\documentclass[11pt]{article}
\usepackage{amscd,amsfonts, amsmath, amsthm, amstext, amssymb, mathrsfs,graphicx,verbatim}

\newtheorem{lemma}{Lemma}[section]
\newtheorem{theorem}[lemma]{Theorem}
\newtheorem{definition}[lemma]{Definition}
\newtheorem{proposition}[lemma]{Proposition}

\newtheorem{remark}[lemma]{Remark}

\title{Isospectral Deformations of Eguchi-Hanson Spaces as   Nonunital Spectral Triples}
\author{\textbf{C. Yang\footnote{E-mail address: chen.yang2@durham.ac.uk}}\\
\textit{}\\
\textit{\normalsize{Department of Mathematical Sciences}}\\
\textit{\normalsize{University of Durham, England, DH1 3LE}}}
 
\begin{document}         
                
\maketitle

\begin{abstract}
We  study the isospectral deformations of the Eguchi-Hanson spaces along a torus isometric action  in the noncompact noncommutative geometry. We concentrate on locality, smoothness and summability conditions of  the nonunital spectral triples, and relate them to the geometric conditions to be   noncommutative spin manifolds.\end{abstract}

\tableofcontents

\def\Xint#1{\mathchoice 
{\XXint\displaystyle\textstyle{#1}}%
{\XXint\textstyle\scriptstyle{#1}}%
{\XXint\scriptstyle\scriptscriptstyle{#1}}%
{\XXint\scriptscriptstyle\scriptscriptstyle{#1}}%
\!\int} 
\def\XXint#1#2#3{{\setbox0=\hbox{$#1{#2#3}{\int}$} 
\vcenter{\hbox{$#2#3$}}\kern-.5\wd0}} 
\def\ddashint{\Xint=} 
\def\dashint{\Xint-}

\section{Introduction}
As a generalization of  Connes' noncommutative differential geometry \cite{Connes-1994}, noncompact noncommutative geometry  is the study of nonunital spectral triples \cite{rennie-2003},\cite{Gayral-Bondia-Iochum-Varilly-2004}.  Various authors also consider the aspect of  summability  as in \cite{estrada-1998},\cite{rennie-2004},\cite{gayral-2006-237}. 

In the unital case, 
Connes provides a set of axioms for  unital spectral triples so to define compact noncommutative spin manifolds. See for example \cite{Connes-1994}.  Rennie and V\'arilly explicitly reconstruct  compact commutative spin manifolds  from  slightly modified axioms \cite{rennie-2006}.
As to the nonunital case, a complete generalization considering these axioms  is not known yet. There are various nonunital examples \cite{graciabondia-2002-0204},\cite{Gayral-Bondia-Iochum-Varilly-2004},\cite{gayral-2007}, which may serve the purpose of testing the axioms or  geometric conditions suggested. In this article, we obtain another nonunital example   by isospectral deformation of  Eguchi-Hanson (EH-) spaces  \cite{eguchi-hanson-1978}. They are   geodesically complete Riemannian spin manifolds in the commutative geometry.

Isospectral deformation is a simple method to deform a commutative spectral triple. It traces back to the Moyal type of deformation from quantum mechanics. Rieffel's insight is to  consider Lie group  actions  on function spaces and hence explain the Moyal product between functions by oscillatory integrals over the group actions \cite{Rieffel-1993}.  Apart from the well-known Moyal planes and  noncommutative tori \cite{rieffel-1990}, this scheme allows more general deformations. Connes and Landi in \cite{Connes-landi-2001} deform spheres and more general compact spin manifold with isometry group containing a two-torus. Connes and Dubois-Violette in \cite{Connes-Dubois-Violette2002} observe that this  works equally well  for noncompact spin manifolds. As in the appendix of \cite{rennie-2003}, it is possible to fit such noncompact examples in the nonunital framework there. 
 The deformation of   EH-spaces we will consider in the following is obtained by these methods and serves as an example of a nonunital triple.
  
The Eguchi-Hanson spaces are of interest  in both Riemannian geometry and physics. Geometrically, they are the simplest asymptotic locally Euclidean  (ALE) spaces, for which a complete classification is provided by  Kronheimer through the method of hyper-K\"ahler quotients \cite{Kronheimer-1989-1}. This construction realizes the family of EH-spaces as 
 a resolution of a singular conifold. In physics, where they first appeared, EH-spaces are known as gravitational instantons. Due to their  hyper-K\"ahler structure, the ADHM construction \cite{Atiyah-1978}, obtaining Yang-Mills' instantons,  is generalized on the EH-spaces in an elegant  way \cite{Kronheimer-1990},\cite{Nakajima-1994}. 
The nonunital spectral triple from isospectral deformation of Eguchi-Hanson spaces may thus link various perspectives.  

Our aim in this article is to concentrate on the locality, smoothness \cite{rennie-2003} and summability conditions of these triples and further see how they fit  into the modified geometric conditions  for nonunital spectral triples. 

The organization of the rest of the article is as follows.  In section 2, we describe the Eguchi-Hanson spaces in the spin geometry. In section 3, we consider  algebras of functions over EH-spaces, the deformation quantization of algebras, and  representations of algebras as operators on the Hilbert space of spinors. We also obtain a projective module description of the spinor bundle. In section 4, we define spectral triples of the deformed EH-spaces and study their summability. In section 5, we discuss how the triple fits into the modified geometric conditions. We conclude in Section 6.


\section{Spin geometry of Eguchi-Hanson spaces}
In this section, we first describe
 the metric and the Levi-Civita connection of the Eguchi-Hanson space,
 and then introduce its spinor bundle, the spin connection and the  Dirac operator. Finally, we  write down the torus action through parallel propagators on the spinor bundle.

\subsection{Metrics, connections and torus isometric actions}

The Eguchi-Hanson spaces were  originally  constructed as gravitational instantons 
 \cite{eguchi-hanson-1978}.   Generalized by Gibbons and Hawking, they fall into a new  category of solutions of  the Einstein's equation,  known as the multicenter solutions \cite{gibbons-hawking-1978}. In local coordinates, the metric is
\begin{equation}
\label{Haw}
ds^2=\Delta
^{-1} dr^2+r^2\left[(\sigma_x^2+\sigma_y^2)+ \Delta
\,  \sigma_z^2\right],
\end{equation}
where $\Delta:=\Delta(r):
=1-a^4/r^4$ and $\{\sigma_x, \sigma_y,\sigma_z\}$ are the standard Cartan basis for three sphere,
\begin{eqnarray*}
\sigma_x&=&\frac{1}{2}\, (-\cos\psi\, d\theta-\sin\theta\, \sin\psi\, d\phi),\\
\sigma_y&=&\frac{1}{2}\, (\sin\psi\, d\theta-\sin\theta\, \cos\psi\, d\phi),\\
\sigma_z&=&\frac{1}{2}\, (-d\psi-\cos\theta\, d\phi),
\end{eqnarray*}
with
 $r\geq a, \quad0\leq\theta\leq\pi, \quad 0\leq\phi< 2\pi, \quad 0\leq\psi< 2\pi.$
 
\begin{remark}
The convention that the period of $\psi$ is $2\pi$ rather than $4\pi$ as in the original construction
is suggested in \cite{gibbons-hawking-1978} to remove the singularity at $r=a$, so that  the  manifold becomes geodesically complete. 
\end{remark}

The  EH-space is diffeomorphic to the tangent bundle of a $2$-sphere $T(\mathbb S^2)$. Modulo a distortion of the  metric, the base as a unit two sphere $\mathbb S^2$ is parametrized by parameters $\phi$ and $\theta$, with $\theta=0$ as the south pole and  $\theta=\pi$ as the north pole. The angle $\phi$ parametrises  the circle defined by a constant $\theta$. Over each point, say $(\theta, \phi)$  on the $2$-sphere, the tangent plane is parametrized by $(r, \psi)$. $r$ parametrizes the radial direction with $r=a$ at the origin of the plane. Circles of constant $r$ are parametrised by $\psi$. 
The identification of $\psi=\psi+2\pi$ is the identification the antipodal points on the circle of constant radius. Together with the metric, this implies that the space at  large enough $r$ is asymptotic to $\mathbb R^4/\mathbb Z^2$, so that it is an ALE space. 

The parameter $a$ in the metric  (\ref{Haw}) is a non-negative real number parametrizing a family of EH-spaces. When $a=0$, the metric degenerates to the conifold $\mathbb R^4/\mathbb Z^2$ and the rest of the family is a resolution of the conifold. This appears  as the simplest case  in Kronheimer's classification of ALE spaces \cite{Kronheimer-1989-1}.  We will  only concentrate on the smooth case so that $a$ is assumed to be positive.

\label{lcc}
Choose the local coordinates $\{x_i\}$ with 
$x_1=r,$ $x_2=\theta,$ $x_3=\phi,$ $x_4=\psi.$
We will write the coordinates $(r, \theta, \phi, \psi)$ and $(x_1, x_2, x_3, x_4)$ interchangeably throughout the article, because the former give a clear geometric picture while the latter are convenient  in tensorial expressions. 
The corresponding basis on the tangent space $T_x(EH)$ of any point $x\in EH$ are  $\left\{\partial_i:=\frac{\partial}{\partial x_i}\right\},$  and  the dual basis on the cotangent space $T_x^*(EH)$ are  
$\{dx^j\}.$
The corresponding metric tensor $g_{ij}(x)\, dx^i\otimes dx^j$  can be written as entries of the matrix $G=(g_{ij})$ as
\begin{equation}
\label{matrixG}
G(x)=\frac{1}{4}\begin{pmatrix}
4\, \Delta
^{-1}& 0                  & 0                                                           & 0\\
0                  & r^2 & 0                                                            &0\\
0                  & 0                  &  \rho
 & r^2\, \Delta
\, \cos\theta  \\
0                   &0                  & r^2\, \Delta
\, \cos\theta    &r^2\, \Delta
 
\end{pmatrix}
\end{equation}
where $\rho:=\rho(r, \theta):
=\left(r^4-a^4 \cos^2\theta\right)/r^2.$  We always assume Einstein's summation convention.

In  the same coordinate chart, the Christoffel symbols of the Levi-Civita connection of (\ref{Haw})
, defined by $\nabla_{i}\partial_j=\Gamma_{ij}^k\partial_k,$
are explicitly,
\begin{eqnarray}
\label{levi}
&&\Gamma^1_{11}=-\frac{\Delta
'}{\Delta
} , \quad\Gamma^1_{22} =-\frac{r\, \Delta}{4}
, \quad\Gamma^1_{33}=-\frac{\Delta\, \rho^+
}{4\, r}
,\nonumber\\
&&
\Gamma^1_{34} =-\frac{r\, \Delta^+
 \Delta
\,\cos\theta}{4} , \quad
\Gamma^1_{44}=-\frac{r\, \Delta^+
 \Delta}{4},
\nonumber \\
&&\Gamma^2_{12}= \frac{1}{r}, \quad\Gamma^2_{33}=-\frac{a^4\sin 2\theta}{2\, r^4}, \quad\Gamma^2_{34}=\frac{\Delta
\, \sin\theta}{2}, \nonumber\\
&&
\Gamma^3_{13}=\frac{1}{r}, \quad\Gamma^3_{23}=\frac{\cot\theta\, \Delta^+
}{2}, \quad\Gamma^3_{24}=-\frac{ \Delta}{2\, \sin\theta},
\nonumber\\
&&
\Gamma^4_{13}= \frac{2\, a^4 \, \cos\theta}{r(r^4-a^4)}, \quad
\Gamma^4_{14}= \frac{\Delta^+}{r\, \Delta}, \quad\Gamma^4_{23}=
=-\frac{ \rho^+
}{2\,r^2\, \sin\theta}, \quad\Gamma^4_{24}=\frac{\cot\theta\,\Delta}{2}. 
\end{eqnarray}
where $$\Delta^+:=\Delta^+(r)
:=1+\frac{a^4}{r^4},\quad \Delta
':=\frac{\partial\Delta}
{\partial r},\quad \rho^+:=\rho^+(r, \theta)
:=\frac{r^4+a^4\, \cos^2\theta}{ r^2}.$$
The identity $\Gamma^{i}_{jk}=\Gamma^{i}_{kj}$,
implied by the torsion free property of the connection, generates another set of symbols and  all the rest of the Christoffel symbols vanish.

The isometry group of the metric (\ref{Haw}) is $(U(1)\times SU(2))/\mathbb Z^2.$ The Killing vector $\partial_{\psi}$ generates the group
 $U(1)/\mathbb Z^2$. 
Another  Killing vector  is $\partial_{\phi}$. Its action on the restriction of the space at $r=a$ is  analogous to  one of the three typical generators of the Lie algebra of the Lie group $SU(2)$   on a standard two-sphere.  
These are the two Killing vectors which define a torus action $\sigma$ on the Eguchi-Hanson space,
\begin{equation}
\label{sigmaact}
\sigma:U(1)\times U(1)\longrightarrow Aut(EH),
\end{equation}
by
$\sigma(\exp{(i\, t_3\,\partial_{\phi})}, \exp{(i\, t_4\, \partial_{\psi})})(r, \theta, \phi, \psi)=(r, \theta, \phi+t_3, \psi+t_4),$
where  $0\leq t_3<2\pi$, $0 \leq t_4<2\pi$ and for any point $ (r, \theta, \phi, \psi)\in EH$. The isometric torus action will determine the   isospectral deformation later.

\subsection{The stereographic projection and orthonormal basis}
\label{stereo}
We choose  an orthonormal basis  to trivialize the cotangent bundle of the EH-space and obtain the corresponding transition functions. 
Since the EH-space is topologically the same as $T(\mathbb S^2)$, we may obtain another set of coordinates by taking the  stereographic projection of the  $\mathbb S^2$ part, while keeping the coordinates on the  tangent space unchanged.

The EH-space (\ref{Haw}) can be covered by two open neighbourhoods $U_N$ and $U_S$, where $U_N$ covers the whole space except  at $\theta=\pi$  and $U_S$ covers the whole space except  at $\theta=0$. 
We may define the map
$f_N:U_N\longrightarrow\mathbb C\times\mathbb R^2$
 by taking a stereographic projection of the base two sphere to $\mathbb C$. I.e., 
 $f_N( \phi, \theta, r, \psi)=(z; r, \psi)$. For the coordinate chart $U_S$, we similarly define the projection map 
$f_S:N_S\longrightarrow\mathbb C\times\mathbb R^2,$
by
 $f_S(\phi, \theta, r, \psi)=(w; r, \psi),$
 where $$z:=\cot\frac{\theta}{2} \, e^{-i\phi},\quad w:=\tan\frac{\theta}{2}\, e^{i\phi}.$$

For any point $ x\in U_N\cap U_S$, the transition function from the coordinate charts $U_S$ to $U_N$ is 
$$(w, r, \psi)=(\frac{1}{z}, r, \psi),$$
and the transition function from $U_N$ to $U_S$ is $(z, r, \psi)=(\frac{1}{w}, r, \psi).$

The  restriction of metric (\ref{Haw}) on the $U_N$ chart with coordinates $(z; r, \psi)$ is,
\begin{eqnarray*}
ds^2&=&\frac{r^2}{(1+z\overline z)^2} dz d\overline z
+\frac{r^2\, \Delta}{4}\left[d\psi+\frac{z\overline z-1}{z\overline z+1}\frac{i}{2}(\frac{dz}{z}-\frac{d\overline z}{\overline z})\right]^2.
\end{eqnarray*}

To obtain a local orthonormal basis of $T^*(EH)_{U_N}$ we may simply define

$$l:=\frac{r}{\sqrt{2}\, (1+z\overline z)} dz, \quad m:=\frac{ \Delta^{-1/2}}{ \sqrt{2}}
 dr+\frac{ r\, \Delta^{1/2}}{4\, \sqrt{2}}
\left[
      \left(\frac{dz}{z}-\frac{d\overline z}{\overline z}\right)\frac{1-z\overline z}{1+z\overline z}+2i\, d\psi
\right],$$
with their complex conjugates
$\overline l$,
$\overline m$ so that the metric tensor over $U_N$ is 
$ds^2=l\otimes\overline l+\overline l\otimes l+m\otimes \overline m+\overline m\otimes m.$

A  real orthonormal frame $\{\vartheta^{\alpha}\}$ of $T^*(EH)_{U_N}$ is thus defined by
$$\vartheta^1:=\frac{1}{\sqrt{2}}(l+\overline l), \quad\vartheta^2:=-\frac{i}{\sqrt{2}}(l-\overline l),\quad
\vartheta^3:=-\frac{i}{\sqrt{2}}(m-\overline m), \quad\vartheta^4:=\frac{1}{\sqrt{2}}(m+\overline m).$$
such that the metric on $U_N$ is diagonalized as
$ds^2=\delta_{\alpha\beta}\vartheta^\alpha\otimes \vartheta^\beta$.
The coordinate transformations  $\vartheta^{\alpha}=h^{\alpha}_i dx^i$ are determined by the  matrix $H=(h^{\alpha}_i)$,
\begin{equation}
\label{matrixH}
H=
\frac{1}{2}\begin{pmatrix}
0                      & -r\, \cos\phi  & -r\, \sin\theta\, \sin\phi                 & 0\\
0                      & r\, \sin\phi    & -r\,\sin\theta\, \cos\phi               & 0\\
0                      & 0                                & r\, \Delta
^{1/2}\, \cos\theta   & r\, \Delta
^{1/2}\\
\frac{2}{\Delta
^{1/2}} & 0                               & 0                                                            & 0
\end{pmatrix},
\end{equation}
whose inverse $H^{-1}=(\tilde h_{\beta}^j)$ from $dx^j=\tilde h^j_{\beta}\vartheta^{\beta}$ is
\begin{equation}
\label{matrixH1}
H^{-1}=2
\begin{pmatrix}
0                                                  &   0            & 0                            &\frac{\Delta
^{1/2} }{2} \\
-\frac{\cos\phi }{r} &\frac{ \sin\phi}{r}                     &0                      &0\\
-\frac{ \sin\phi}{r\,\sin\theta}&  -\frac{ \cos\phi  }{r\,\sin\theta}             & 0&0\\
\frac{\cos\theta\, \sin\phi}{r\,\sin\theta}& \frac{\cos\theta\, \cos\phi }{r\,\sin\theta}    &\frac{1}{r\, \Delta
^{1/2}}
&0.
\end{pmatrix}
\end{equation}

The above construction on the $U_N$ chart can be carried out  the same way on the $U_S$ coordinates. We denote orthonormal frames over $U_S$ by adding $'$'s to $l, m, \vartheta^{\alpha}$, $x_j$ and etc. 


Local frames $\{\vartheta^{\alpha}\}$ on $U_N$ define  a local trivialization of the cotangent bundle,
$F_N: T^*(EH)_{U_N}\rightarrow U_N\times\mathbb R^4$
by $F_N(x; a_1 \vartheta^1+\cdots+a_4 \vartheta^4):=(x;  a_1,\dots, a_4)$, 
where $a_{\alpha}$'s  are real-valued functions over $U_N$. 
In a similar way, the choice of local frames $\{\vartheta'^{\alpha}\}$ on $U_S$ defines a local trivialization of the cotangent bundle,
$F_S: T^*(EH)_{U_S}\longrightarrow U_N\times\mathbb R^4$.


The transition functions $f^{\beta}_{\alpha}$'s such that $\vartheta'^{\beta}=f^{\beta}_{\alpha} \vartheta^{\alpha}$ are elements of the matrix
 $F_{SN}:=F_N\circ F_S^{-1}$ as
\begin{equation}
\label{matrix}
F_{SN}=\begin{pmatrix}
-\frac{\overline z^2+z^2}{2\,z\overline z} &  -i\,\frac{\overline z^2-z^2}{2\, z\overline z} & 0&0\\
i\, \frac{\overline z^2-z^2}{2\, z\overline z}& -\frac{\overline z^2+z^2}{2\, z\overline z} & 0 &0 \\
0 &0 &1 &0\\
0 &0 & 0 &1
\end{pmatrix}
=
\begin{pmatrix}
-\cos 2\phi & \sin 2\phi & 0&0\\
-\sin 2\phi &-\cos 2\phi & 0 &0\\
0 &0 &1 &0\\
0 &0 & 0 &1
\end{pmatrix}.
\end{equation}
The inverse  transition function  is given by the inverse of the matrix $F_{SN}$,
$F_{NS}:=F_S\circ F_N^{-1}=F_{SN}^{-1}.$
The cotangent bundle is thus
\begin{equation}
\label{cotang}
T^*(EH)=(U_N\times\mathbb R^4)\cup (U_S\times\mathbb R^4)/\sim,
\end{equation}
where 
 $(x;  a_1, \dots, a_4)\in U_N\times\mathbb R^4$ and  $( x';  a'_1,\dots, a'_4)\in U_S\times\mathbb R^4$ are  defined to be equivalent 
if and only if $ x= x'$ and $F_{NS}(a_1,\dots, a_4)^t=(a'_1, \dots, a'_4)^t$.

\subsection{Spin structures and spinor bundles}
\label{EHS}
Following a standard procedure from \cite{Lawson-1989}, we obtain the spinor bundle of the EH-space.  In  coordinate charts $\{U_N, U_S\}$,
the frame bundle $P_{SO(4)}$ of the EH-space is  the $SO(4)$-principal bundle with transition functions
 $F_{NS}$  in (\ref{matrix}) and its inverse $F_{SN}$. 

Recall that the   covering  map of groups,
\begin{equation}
\label{double}
\rho: Spin(4)\longrightarrow SO(4),
\end{equation}
 is defined by  the adjoint representation of  $Spin(4)$ as
$\rho(w) x:=w\cdot x\cdot w^t$ for $x\in\mathbb R^4$,
where $w=v_1\cdots v_{m}\in Spin(4), $ $m$ is even and  $v_i\in\mathbb R^4$ for $i=1, \dots, m$. Geometrically,  $\rho(w)
=\rho(v_1)\circ\cdots\circ \rho(v_{m}),$ where
$\rho(v_i)$ is the reflection of the space $\mathbb R^4$ with respect to the hyperplane with normal vector  $v_i$. 

Locally, the upper left block of the transition  matrix (\ref{matrix})  is  a rotation in the plane spanned by $\{\vartheta^1, \vartheta^2\}$ through an angle $2\phi+\pi$. Such rotation can be decomposed to  two reflections say $\rho(v_2)\circ\rho(v_1)$, with
$$v_1:=\vartheta^1,\quad v_2:=-\sin\phi\, \vartheta^1+\cos\phi\, \vartheta^2.$$
\begin{remark}
Another choice is $\rho(-v_2)\circ\rho(v_1)$, which gives the same rotation as an element in $SO(2)$. 
\end{remark}

 $v_2\cdot v_1\in Spin(4)$ is a lifting of $\rho(v_2)\circ\rho(v_1)\in SO(4)$ under the covering map (\ref{double}). Thus, in the local coordinate chart $U_N$, $\widetilde{F_{SN}}:=v_2\cdot v_1$ in $Spin(4)$ defines a lifting of the action $F_{NS}\in SO(4)$ as in (\ref{matrix}) under the double covering (\ref{double}).
 
 To obtain a global lifting of the frame bundle, we consistently define the transition matrix $\widetilde{F_{NS}}$ as a lifting in the group $Spin(4)$  over  $ x'\in U_S$ by
$\widetilde{F_{NS}}=-v_2'\cdot v_1',$ where
$$v_1':=\vartheta'^1,\quad v_2':=\sin\phi \, \vartheta'^1+\cos\phi\, \vartheta'^2. $$
The following confirms the consistency  of the liftings on two coordinate charts.
\begin{lemma}
Transition functions $\{\widetilde{F_{NS}}, \widetilde{F_{SN}}\}$ satisfy the cocycle condition, $\widetilde{F_{NS}}\circ\widetilde{F_{SN}}=\widetilde{F_{SN}}\circ\widetilde{F_{NS}}=\mathbf 1.$
\end{lemma}
\begin{proof}
Applying the transformation from $\vartheta^{\alpha}$'s to $\vartheta'^{\beta}$'s  by (\ref{matrix}), we have $\vartheta'^1\cdot\vartheta'^2=\vartheta^1\cdot\vartheta^2$. Thus,
\begin{eqnarray*}
\widetilde{F_{NS}}\circ\widetilde{F_{SN}}&=&-v_2'\cdot v_1'\cdot v_2\cdot v_1\\
&=&-(\sin\phi\, \vartheta'^1+\cos\phi\, \vartheta'^2)\cdot\vartheta'^1\cdot(-\sin\phi\, \vartheta^1+\cos\phi\, \vartheta^2)\cdot\vartheta^1\\
&=&\sin^2\phi-\sin\phi\,\cos\phi(\vartheta^1\cdot\vartheta^2+\vartheta^2\cdot\vartheta^1)-\cos^2\phi\, \vartheta^2\cdot\vartheta^1\cdot \vartheta^2\cdot\vartheta^1=\mathbf 1,
\end{eqnarray*}
by using identities $\vartheta^{\alpha}\cdot\vartheta^{\alpha}=\mathbf 1$ and $\vartheta^{\alpha}\cdot\vartheta^{\beta}=\vartheta^{\beta}\cdot\vartheta^{\alpha}$ for $\alpha\neq\beta$, of elements of the orthonormal bases $\vartheta^{\alpha}$'s and those of $\vartheta'^{\beta}$'s. Similarly, $\widetilde{F_{SN}}\circ \widetilde{F_{NS}}=\mathbf 1$.

\end{proof}
Therefore,
the principal $Spin(4)$-bundle can be  defined by
\begin{equation}
\label{spin4}
P_{Spin(4)}:=(U_N\times Spin(4)\cup U_S\times Spin(4))/\sim.
\end{equation}
where  
$(x, \tilde g)\in U_N\times Spin(4)$ and $(x', \tilde g')\in U_S\times Spin(4)$ are defined to be equivalent if
 and only if $x=x'$ and $\tilde g'=\tilde F_{NS} \, \tilde g$.

The  double covering   of bundles  (\ref{spin4}) over  the EH-space defines a \textit{spin structure} of it. We will always assume this choice of spin structure. 

The spinor bundle can be defined
 as an associative bundle of typical fiber $\mathbb C^4$ of the principal $Spin(4)$-bundle (\ref{spin4}), by  specifying a representation of $Spin(4)$ on $GL_{\mathbb C}(4)$.
We know that locally, for any $ x\in EH$,  there exists a unique irreducible  representation space $\Lambda$ of complex dimension $4$ of the Clifford algebra  $Cl(T_x^*(EH))$ through the Clifford action $c:Cl(T_x^*(EH))\rightarrow End(\Lambda)$.
We define the representation of $Spin(4)$ in $End(\Lambda)(\cong GL_{\mathbb C}(4))$ simply by the restriction of $c$ from the Clifford algebra, and obtain the \textit{spinor bundle} $\mathcal S$ of typical fiber $\Lambda$, with transition functions $\{c(\widetilde{{F}_{NS}}), c(\widetilde{F_{SN}})\}$ in the coordinate charts $\{U_N, U_S\}$.

With respect to the orthonormal basis, say $\{\vartheta^{\alpha}\}$ of $T^*(EH)_{U_N}$, 
there exists a unitary frame $\{f_{\alpha}\}$ of the representation space $\Lambda\cong \mathbb C^4$, 
such that the Clifford representations
$\gamma^{\alpha}:=c(\vartheta^{\alpha}( x))$ for $\alpha=1, \dots, 4$ can be represented as constant matrices,
$$
 \gamma^1=\begin{pmatrix}
0  &  0  &  -1&  0\\
 0  &   0   & 0    &  1 \\
   1&   0  &  0    &  0   \\
  0&   -1   &  0    &  0   \\
   \end{pmatrix},\quad 
   \gamma^2=\begin{pmatrix}
0  &  0  &  -i &  0\\
 0  &   0   & 0 &  -i  \\
-i &   0   &   0   &    0 \\
 0  &    -i  &  0    &   0  \\
   \end{pmatrix},$$
   \begin{equation}
\label{pauli}
 \gamma^3=\begin{pmatrix}
0  &  0  &  0 &  -1\\
  0 &   0   &  -1    &  0  \\
  0 &   1   &   0   &0     \\
  1 &   0   &     0 &  0   \\
   \end{pmatrix},\quad
   \gamma^4=\begin{pmatrix}
0  &  0  &  0 &  -i\\
 0  &   0   &  i  &    0 \\
  0  &  i   &   0   &   0  \\
  -i &   0   &   0   &   0  \\
   \end{pmatrix}.  \end{equation}
The fact is that there exists frames $\{f'_{\beta}\}$ on the coordinate chart $U_S$ so that  the representation of $c(\vartheta'^{\beta})$ are also constant matrices $\gamma^{\beta}$'s as above.



Under the chosen frames $\{f_{\alpha}\}$ and $\{f'_{\beta}\}$, we may represent the transition functions of the spinor bundle as follows. 
Define  maps
$P, Q: U_N\cap U_S\longrightarrow GL_{\mathbb C}(4)$
 by
\begin{equation}
\label{qsn}
P:
=c( \widetilde{F_{SN}})
=-\sin\phi\gamma^1\gamma^1+\cos\phi\gamma^2\gamma^1
= diag(- \frac{i\,\overline z}{|z|},  \frac{i\,z}{|z|},  \frac{i\,z}{|z|}, - \frac{i\,\overline z}{|z|}),
\end{equation}
\begin{equation}
\label{qns}
Q:
=c( \widetilde{F_{NS}})
=-\sin\phi\gamma^1\gamma^1-\cos\phi\gamma^2\gamma^1
=diag(\frac{i\,\overline w}{|w|}, -\frac{i\, w}{|w|}, -\frac{i\, w}{|w|}, \frac{i\,\overline w}{|w|}).
\end{equation}
$diag(a, b, c, d )$ stands for diagonal matrix with diagonal elements $a, b, c, d$.

The 
spinor bundle  $\mathcal S$ is thus,
\begin{equation}
\label{spinorbdle}
\mathcal S:=(U_N\times \mathbb C^4\cup U_S\times\mathbb C^4)/\sim,
\end{equation}
where   $( x; s_1,\cdots, s_4)\in U_N\times \mathbb C^4$ and  $( x'; s_1', \cdots, s_4')\in U_S\times \mathbb C^4$ are defined to be equivalent 
if and only if $ x= x' $ and $(s_1', \cdots, s_4')^t= Q
(s_1, \cdots, s_4)^t$. One can easily see that the cocycle condition of the transition functions $P
\circ Q
=Q
\circ P
=\mathbf 1$ holds.


The \textit{chirality operator} is defined
by
   \begin{equation}
   \label{chi}
   \chi:=c(\vartheta^1)\,  c(\vartheta^2)\,  c(\vartheta^3)\,  c(\vartheta^4)=\gamma^1\, \gamma^2\, \gamma^3\, \gamma^4=diag(-1, -1, 1, 1),
   \end{equation}
 such that $\chi^2=\mathbf 1$. The representation space
   $
   \Lambda=\Lambda^+\oplus\Lambda^-
 $
   is decomposed as $\pm 1$-eigenspaces of  the operator $\chi$, with
   $dim_{\mathbb C}\, \Lambda^+=dim_{\mathbb C}\, \Lambda^-=2$. 
   This fiberwise splitting extends to  the global decomposition of spinor bundle as subbundles over the EH-space,
$
\mathcal S=\mathcal S^+\oplus\mathcal S^-,
$
with each of the complex subbundles $\mathcal S^+$ and $\mathcal S^-$  of rank $2$.
Therefore,  any element $s\in\mathcal S$ can be decomposed as
$s=(s^+, s^-)^t$.
The \textit{charge conjugate operator} on the spinor bundle $J:\mathcal S\rightarrow \mathcal S$ is defined by
\begin{equation}
\label{charge}
J
\begin{pmatrix}
s^+\\
s^-
\end{pmatrix}
:=
\begin{pmatrix}
-\overline{s}^-\\
\overline{s}^+
\end{pmatrix}.
\end{equation}


\subsection{Spin connections and Dirac operators of  spinor bundles}
\label{sectionspin}
Following the general procedure in \cite{Gracia-Bondia-2001}, we can induce the spin connection $\nabla^{\mathcal S}$ of the spinor bundle $\mathcal S$ from the Levi-Civita connection of the EH-space. 

We will only work on the $U_N$ coordinate chart and the construction  on $U_S$ is similar. 
In the orthonormal frame $\{\vartheta^{\alpha}\}$,  the corresponding Levi-Civita connection on the dual tangent bundle, $T^*(EH)_{U_N}$,  can be expressed as
$
\nabla^{T^*EH}\vartheta^{\beta}=-\widetilde{\Gamma}_{i\alpha}^{\beta}\, dx^i\otimes\vartheta^{\alpha}.
$
The metric compatibility of the Levi-Civita connection implies that $\widetilde{\Gamma}^{\alpha}_{i\beta}=-\widetilde{\Gamma}_{i\alpha}^\beta$.

We may represent $\widetilde{\Gamma}_{i\alpha}^{\beta}$'s in terms of the Christoffel symbols $\Gamma_{ij}^k$'s of $\nabla$ in the $dx^i$'s (\ref{levi}) by
\begin{equation}
\label{Chris}
\widetilde{\Gamma}_{i\alpha}^\beta=\tilde h^{j}_{\alpha}\left(h_k^{\beta}\, \Gamma_{ij}^k-\partial_i h_j^{\beta}\right),
\end{equation}
where $h_i^{\alpha}$'s and  $\tilde h_{\beta}^j$'s are the matrix entries of $H$ in (\ref{matrixH}) and  $H^{-1}$ in (\ref{matrixH1}), respectively. 
Modulo the anti-symmetric condition between $\alpha$ and $\beta$ indices, all  the nonvanishing Christoffel symbols are
\begin{eqnarray}
\label{Gamma}
&&\widetilde{\Gamma}_{23}^{1}=\frac{1}{2}\, \Delta
^{1/2}\, \sin\phi, \quad\widetilde{\Gamma}_{24}^1=-\frac{1}{2}\, \Delta
^{1/2}\, \cos\phi,
\nonumber\\
&&\widetilde{\Gamma}_{22}^{3}=\frac{1}{2}\, \Delta
^{1/2}\, \cos\phi, \quad
\widetilde{\Gamma}_{22}^4=-\frac{1}{2}\, \Delta
^{1/2}\, \sin\phi, \nonumber\\
&&\widetilde{\Gamma}_{33}^{1}=-\frac{1}{2}\, \Delta
^{1/2}\, \sin\theta\, \cos\phi, \quad \widetilde{\Gamma}_{34}^{1}=-\frac{1}{2}\, \Delta
^{1/2}\, \sin\theta\, \sin\phi,\nonumber \\
&&
\widetilde{\Gamma}_{32}^{1}=-1-\frac{1}{2}\, \Delta^+\, \cos\theta, \quad
\widetilde{\Gamma}_{32}^{3}=-\frac{1}{2}\, \Delta
^{1/2}\, \sin\theta\, \sin\phi\nonumber\\
&&
\widetilde{\Gamma}_{32}^{4}=\frac{1}{2}\, \Delta
^{1/2}\, \sin\theta\, \cos\phi, \quad
\widetilde{\Gamma}_{33}^{4}=-\frac{1}{2}\, \Delta^+\, \cos\theta\nonumber\\
&&
\widetilde{\Gamma}_{42}^{1}=\frac{1}{2}\, \Delta
, \quad\widetilde{\Gamma}_{43}^{4}=-\frac{1}{2}\, \Delta^+.
\end{eqnarray}

We  define $\gamma_{\alpha}:=\gamma^{\alpha}$, then
the \textit{spin connection} $\nabla^{\mathcal S}:\mathcal S\rightarrow \mathcal S\otimes \Omega^1(EH)$ is 
\begin{equation}
\label{spinconn}
\nabla^{\mathcal S}:=d-\frac{1}{4}\widetilde{\Gamma}_{i\alpha}^{\beta}\, dx^i\otimes\gamma^{\alpha}\gamma_{\beta}.
\end{equation}



The covariant derivative   $\nabla^{\mathcal S}_{ i}:=\nabla^{\mathcal S}(\partial_i)$,   for $i=1, \dots, 4$, equals to 
$\nabla_{i}^{\mathcal S}=\partial_i-\omega_i, $ where $\omega_i=\frac{1}{4}\widetilde{\Gamma}^{\beta}_{i\alpha}\gamma^{\alpha}\gamma_{\beta}.$
The \textit{Dirac operator}  $\mathcal D:\Gamma(\mathcal S)\rightarrow\Gamma(\mathcal S)$ 
can be defined by
\begin{equation}
\label{Diracx}
\mathcal D(\psi):=-i\,  \gamma^j \, \nabla_{j}^{\mathcal S}\psi,\quad\forall\psi\in \Gamma(\mathcal S),
\end{equation}
where $\gamma^j:=c(dx^j)=\tilde h^j_{\beta}\gamma^{\beta}$. We note that the compatibility of the spin connection with respect to the spin structure implies that the commutativity between the Dirac operator and the charge conjugate operator, i.e. $[\mathcal D, J]=0$. 



\subsection{Torus actions on the spinor bundle}
A torus action on the spinor bundle $\mathcal S$ can be 
induced from the torus isometric action on a general Riemannian manifold \cite{Connes-landi-2001},\cite{Connes-Dubois-Violette2002}. 
In this subsection,  we will represent such torus action  
 (\ref{sigmaact})  through parallel transporting   spinors along  geodesics. 

Recall that the isometric action $\sigma$ is generated by the two Killing vectors $\partial_{3}=\partial_{\phi}$ and $\partial_4=\partial_{\psi}$. 
Let $c_k:\mathbb R\rightarrow EH$   be the geodesics obtained as integral curves of the Killing vector field $\partial_k$ for  $k=3,4$.

The equation of  parallel transportation with respect to the spin connection  along any curve $c(t)$ is $\nabla^{\mathcal S}_{c'(t)}\psi=0,$ where $c'(t):=d c(t)/dt$,  for $\psi\in\Gamma(\mathcal S).$ 
Substituting (\ref{spinconn}), we obtain
\begin{equation}
\label{matrixA} 
\frac{d\psi}{dt}-A(c(t))\,\psi=0, \quad A(c(t)):=\frac{1}{4}\widetilde{\Gamma}_{i\alpha}^{\beta}\, dx^i(c'(t))\otimes\gamma^{\alpha}\gamma_{\beta}.
\end{equation} 

When the curve is $c_3(t)$, the corresponding matrix $A(c_3(t))$ is  
\begin{equation}
\label{A3}A(c_3(t))=\frac{1}{2}
\begin{pmatrix}
i& 0 &0 & 0\\
0 &-i &0 &0\\
0 & 0 &-i\, (1+\Delta^+\cos\theta) &-\Delta^{1/2}\, \sin\theta\, e^{i\phi}\\
0 & 0& \Delta^{1/2}\sin\theta\, e^{-i\phi} &i\, (1+\Delta^+ \cos\theta),
\end{pmatrix}
\end{equation}
where $r, \theta$ and $\phi$  are understood as components of coordinates on the curve $c_3(t)$. 
When the curve is $c_4(t)$, the corresponding matrix $A(c_4(t))$  is 
\begin{equation}
\label{A4}
A(c_4(t))=\frac{i}{2}\, diag(-\frac{a^4}{r^4}, \frac{a^4}{ r^4}, -1, 1),
\end{equation}
where $r$ is understood as one of the components of coordinates on the curve $c_4(t)$.

The corresponding \textit{parallel propagator}  is a map $P_{c(t)}(t_{0}, t_1):\Gamma(\mathcal S)\rightarrow\Gamma(\mathcal S)$ defined by parallel transporting  any section $\psi$  along the curve $c(t)$ with $t\in[t_0, t_1]$. 
The propagator can be represented by an iterated integration of the equation (\ref{matrixA}). For geodesics $c_k(t)$, $k=3, 4$, the corresponding  matrix is formally solved as
\begin{equation}
\label{propagator}
P_{c_k(t)}(t_0, t_1)=\mathcal P\exp\left(\int_{t_i}^{t_f} A_k(t) dt\right), 
\end{equation}
where $\mathcal P$ is the path-ordering operator.

Let $\mathcal H$ be the Hilbert space completion with respect to the $L^2$-inner product on the space of $L^2$-integrable sections of the spinor bundle $\mathcal S$. The parallel propagators  (\ref{propagator}) can be extended to families of  operators $U_k(t_k-t_0): \mathcal H\rightarrow \mathcal H$ parametrized by the real number  $(t_k-t_0)$  by
\begin{eqnarray*}
U_k(t_k-t_0)(\psi)(x):=(P_{c_k(t)}(t_0, t_k)\psi)(x),\quad\forall \psi\in \mathcal H, 
\end{eqnarray*}
where we assume $x=c_k(t_0)\in EH$  for $k=3, 4$. Without loss of generality, we may take $t_0=0$ so that the family of operators is parametrized by $t_k$. 

Since the  spin connection is compatible with the metric of the EH-space,   the pointwise inner product of the images of any two sections under parallel transportation along  the geodesics $c_k(t)$ remains unchanged. This further implies 
 that their  $L^2$-integrations remain the same. Therefore, the operators $U_k(t_k)$ are unitary. 
Let $W_k$  be the self-adjoint operators on $\mathcal H$ which generate  $U_k$  by $U_k(t_k)=\exp(it_k W_k),$ where $ t_k\in \mathbb R$ for  $k=3, 4.$

We may define a representation of the  
 double cover $p:\widetilde{\mathbb T^2}\rightarrow\mathbb T^2$ 
of the two  torus by $\widetilde{V}:\widetilde{\mathbb T^2}\rightarrow \mathcal L(\mathcal H)
$
such that
\begin{equation}
\label{actiondb}
\widetilde{ V}(\tilde t_3, \tilde t_4)\psi(x):=
\exp(i(\tilde t_3 W_3(x)+\tilde t_4 W_4(x))\psi(x), \quad\forall \psi\in \mathcal H.\end{equation}
This action  covers the isometric action  $\sigma$ of $\mathbb T^2$ from (\ref{sigmaact})
 in the sense that for any $\tilde v\in \widetilde{\mathbb T}^2$,  $p(\tilde v)=v$ implies 
  $$\widetilde{V}_{\tilde v}(f\psi)=\alpha_v(f) \widetilde{V}_{\tilde v}(\psi),\quad\forall \psi\in \mathcal H,$$
 for any bounded continuous function $f\in C_b(EH)$ and the action $\alpha$  on $C_b(EH)$  defined by $\alpha_v(f)(x):=f(\sigma_{-v}(x))$. We assume the  choice of the lifting in the double torus is always fixed and  omit the $\tilde{\cdot}$  for notational simplicity  from now on.

\section{Smooth algebras and projective modules}
\label{smoothalg}
We consider  algebras of functions over the Eguchi-Hanson spaces, and their deformations as  differential algebras.  To obtain a $C^*$-norm on the deformed algebra, we consider  representations of algebras as operators on the Hilbert space of spinors. Some algebras may be realized as  smooth algebras \cite{rennie-2003}.  We also find  projective modules from  the spinor bundle. 


\subsection{Algebras of smooth functions}
We first summarize some related facts on topological algebras of complex-valued functions in \cite{rennie-2003}.
For  a noncompact Riemannian manifold $X$,  let $C_c^{\infty}(X)$  be the space of smooth  functions on $X$ of compact support, 
$C_0^{\infty}(X)$  be the space of smooth  functions  vanishing at infinity and
 $C_b^{\infty}(X)$ be the space of smooth functions  whose derivatives are bounded to all degrees.


In some  local coordinate charts with corresponding partition of unity, say $\mathcal U
=\{U_a, h_a\}_{a\in A}$, we may define the  family of seminorms on $C_b^{\infty}(X)$
by 
\begin{equation}
\label{semi1}
q_m^{\mathcal U
}(f):=\sum_{a\in A}\sup_{| \alpha|\leq m}\left(\sup_{x\in U_a}|h_a(x) \, \partial^{\alpha}f(x)|\right).
\end{equation}
for any   $f\in C_b^{\infty}(X)$, $\alpha$ are multi-indices and $m$ a non-negative integer. These seminorms  restrict on $C_0^{\infty}(X)$ and $C_c^{\infty}(X)$. 
The natural topology induced by (\ref{semi1}) is the \textit{topology of uniform convergence of all derivatives}. We can show that two such families of seminorms defined by different coordinate charts are equivalent. Thus the topology defined does not depend on the choice of coordinates $\mathcal U$. We also note that the $q_0$ seminorm in the family of seminorms is nothing but the suprenorm $\|\cdot\|_{\infty}$, which is a $C^*$-norm with the involution defined by normal complex conjugation. 

Algebras $C_b^{\infty}(X)$ and $C_0^{\infty}(X)$ are both Fr\'echet in the topology of uniform convergence of all derivatives, while the algebra $C_c^{\infty}(X)$ is not complete. However, $C_c^{\infty}(X)$ is complete  in the \textit{topology of inductive limit} as the inductive limit of  the   topology obtained by restriction on a family of algebras $C_c^{\infty}(K_n)$, where $\{K_n\}_{n\in \mathbb N}$ is an increasing family of compact subsets in $X$. 

The algebra $C_c^{\infty}(X)$ is dense in the Fr\'echet algebra $C_0^{\infty}(X)$. 

To consider algebras of smooth functions of the Eguchi-Hanson spaces,  we may  use  the coordinate charts $\mathcal U=\{U_N, U_S\}$  defined in Section \ref{stereo}, 
with a partition of unity $\{h_N, h_S\}$ subordinated to them.
The family of seminorms (\ref{semi1}) can be written as
\begin{eqnarray}
\label{Q2}
q_m^{\mathcal U}(f)&=& \sup_{|\alpha|\leq m}\sup_{x \in U_N}|h_N(x)\, \partial^{\alpha} f(x)|+\sup_{|\alpha'|\leq m}\sup_{x' \in U_S}|h_S(x')\, \partial^{\alpha'} f(x')|.
\end{eqnarray}
We obtain the corresponding topological algebras  by taking $X=EH$.

\subsection{Algebras of integrable functions}
\label{integrable}
Apart from  algebras of functions which can be represented as operators,  there are algebras of functions which  may define projective modules as representation spaces. Decay conditions  at infinity and integrability conditions of  functions become important when considering noncompact spaces. We consider the following algebras of integrable  functions.

The \textit{$(k, p)$-th Sobolev norm} of a function $f$, say in $ C_b^{\infty}(EH)$,  is given as
\begin{equation}
\label{sob}
\|f\|_{ H_k^p}:=\sum_{m=0}^k\left( \int_{EH} |\nabla^m f|^p dVol\right)^{1/p}, 
\end{equation}
where $k$ is a non-negative integer and $p$ is a positive integer. (We will not consider the case where $p$ is a real number).  We define  subspaces 
in $C_b^{\infty}(EH)$ which contain functions with finite Sobolev norm,
$$ C_k^p(EH):=\{f\in C_b^{\infty}(EH): \|f\|_{ H_k^p}<\infty\}.$$
Let $ H_k^p(EH)$ the Banach space obtained by the  completion of the algebra $C_k^p(EH)$ with respect to the Sobolev norm. 
In particular, $ H_0^p(EH)\supset\dots\supset  H_{k}^p(EH)\supset H_{k+1}^p(EH)\supset\cdots.$

\begin{remark}
Notice that the algebra $C_c^{\infty}(EH)$ is contained in  $ H_k^p(EH)$ for any $k\in\mathbb N$. Completion of $C_c^{\infty}(EH)$ with respect to  $\|\cdot\|_{H_k^p}$  gives us 
the Banach space, $ H_{k, 0}^p(EH)$ such that $ H_{k, 0}^p(EH)\subset  H_k^p(EH).$
The equality does not hold in general. However, in the circumstances of a complete Riemannian manifold with Ricci curvature bounded up to degree $k-2$, and positive injective radius (which is  satisfied by the $EH$-space), $H_{k, 0}^p(EH)= H_{k}^p(EH)$ when  $k\geq 2$ \cite{hebey-1999}.
\end{remark}

\begin{lemma}
\label{complete}
For a fixed non-negative integer $p$, the intersection defined as
$$C^{\infty}_p(EH):=\cap_{k} H_k^p(EH)$$
 is a Fr\'echet algebra in the topology defined by the family of norms $\{\|\cdot\|_{ H_k^p}\}_{ k\in\mathbb N}$. 
\end{lemma}
\begin{proof}
The topology is easily seen to be locally convex and metrizable. To show  that it is complete, let $\{f_{\beta}\}$  be 
 any Cauchy sequence in $ C^{\infty}_p(EH)$, then there exists a limit $f_{k}^p$ of $\{f_{\beta}\}$ under the norm $\|\cdot\|_{H_k^p}$ in $ H_k^p(EH)$ for each $k\in\mathbb N$. For any two indices $k_1, k_2$ such that $k_1\leq k_2$, the norm $\|\cdot\|_{H_{k_2}^ p}$ is stronger than the norm $\|\cdot\|_{ H_{k_1}^p}$. The Cauchy sequence $\{f_{\beta}\}$ with the limit $f_{k_2}^p$ in the norm $\|\cdot\|_{ H_{k_2}^p}$ is also  a Cauchy sequence with the limit
$f_{k_1}^p$ in the norm $\|\cdot\|_{ H_{k_1}^p}$. Uniqueness of the limit implies that $f_{k_2}^p=f_{k_1}^p$. Since $k_1, k_2$ are arbitrary,   the limits $f_k^p$  for any $k\in\mathbb N$  agree. We denote the limit as $f$ so that the Cauchy sequence converges to $f\in C_p^{\infty}(EH)$ with respect to any of the norms. Thus the topology is complete and  $C_p^{\infty}(EH)$ is a Fr\'echet algebra.
\end{proof}

When $p=2$,  the Fr\'echet algebra $C^{\infty}_2(EH)$ belongs to the chain of
  continuous inclusions,
\begin{equation}
\label{algebra}
C_c^{\infty}(EH)\hookrightarrow C^{\infty}_2(EH)\hookrightarrow C_0^{\infty}(EH),
\end{equation}
with respect to their aforementioned topologies.

\subsection{Deformation quantizations of differentiable Fr\'echet  algebras}
Rieffel's deformation quantization of a differentiable Fr\'echet algebra in \cite{Rieffel-1993} (Chapter 1, 2) can be summarized as follows.
Let $\mathcal A$ be a Fr\'echet algebra whose topology is defined by a family of seminorms $\{q_m\}$.  We assume that there there is an isometric action $\alpha$ of the vector space $V:=\mathbb R^d$ considered as a $d$-dimensional  Lie group  acting on $\mathcal A$.
We  also assume that the algebra   is smooth  with respect to the action $\alpha$,  i.e. $\mathcal A=\mathcal A^{\infty}$ in the notation of the reference.

Under the choice of a basis $\{X_1, \dots, X_d\}$ of the Lie algebra of $V$, the action  $\alpha_{X_i}$ of $X_i$ defines a partial differentiation on $\mathcal A$.
One can define a new family of seminorms from $q_m$ by taking into account the action of $\alpha$.  For any $f\in \mathcal A$,
\begin{equation}
\label{Q31}
\|f\|_{j, k}:=\sum_{m\leq j,\, |\mu|\leq k} q_m(\delta^{\mu}f),
\end{equation}
where $\mu$ are  the multi-indices $(\mu_1, \dots, \mu_d)$
and $\delta^{\mu}=\alpha_{X_{1}}^{\mu_1}\dots\alpha_{X_d}^{\mu_d}$. The deformation quantization of the  algebra $\mathcal A$ can be carried out in three steps:

\textit{ Step 1. 
}
Let $C_b(V\times V, \mathcal A)$ to be the space of bounded continuous functions from $V\times V$ to $\mathcal A$. One can induce the family of seminorms $\{\|\cdot\|^{C}_{j, k}\}$ on the space $C_b(V\times V, \mathcal A)$ by 
\begin{equation}
\label{Q4}
\|F\|^{C}_{j, k}:=\sup_{w\in V\times V}\|F(w)\|_{j, k},
\end{equation}
for $F$ in $C_b(V\times V, \mathcal A)$ and $\|\cdot\|_{j, k}$ on $\mathcal A$  as in (\ref{Q31}). 

Let $\tau$ be an action of $V\times V$ on the space $C_b(V\times V, \mathcal A)$ defined by translation. That is, $\tau_{w_0}(F)(w)=F(w+w_0)$ for any $w_0, w\in V\times V$ and  $F\in C_b(V\times V, \mathcal A)$.  The action $\tau$ is an isometry action with respect to the seminorms  (\ref{Q4}). We define  $\mathcal B^{\mathcal A}(V\times V)$ to be the maximal subalgebra such that  $\tau$ is  strongly continuous and whose elements  are all smooth with respect to the action of $\tau$.

In the same way as one induces from the family of seminorms $\{q_m\}$
and obtains  the seminorms $\|\cdot\|_{j, k}$ of $\mathcal A$  in (\ref{Q31}), one may induce the family of seminorms on $\mathcal B^{\mathcal A}(V\times V)$ from (\ref{Q4}) by taking into account of the action of $\tau$. For any $F\in \mathcal B^{\mathcal A}(V\times V)$, let
\begin{equation}
\label{qq5}
\|F\|^{\mathcal B}_{j, k; l}:=\sum_{(l, m)\leq(j, k)}\sum_{|\nu|\leq l}\|\delta^{\nu} F\|_{l, m}^C,
\end{equation}
where $\nu$ are the multi-indices and $\delta^{\nu}$ denotes the partial differentiation operator associated to $\tau$ of $V\times V$. 

\textit{Step 2. 
}
The following is the fundamental result of the deformation quantization of a differentiable algebra. See Proposition 1.6 in \cite{Rieffel-1993}. 
For any  invertible map $J$ on $V$, one can define an $\mathcal A$-valued 
oscillatory integral over $V\times V$
of  $F\in \mathcal B^{\mathcal A}(V\times V)$ by \begin{equation}
\label{integration}
\int_{V\times V}F(u, v)  e(u\cdot v) \,du dv,
\end{equation}
where $e(t):=\exp(2\pi i\, t)$ for  $t\in\mathbb R$ and $u\cdot v$ is the natural inner product on $V$ considered as its own Lie algebra.

It is shown to be  $\mathcal A$-valued by getting the bound of the integral in the family of seminorms $\{\|\cdot\|_{j, k}\}$ on $\mathcal A$. Specifically, for large enough $l$,  there exists a constant $C_l$ such that
$$\left\|\int_{V\times V} F(u, v) e(u, Jv)\, du dv\right\|_{j, k}\leq C_l\, \|F\|_{j, k; l}^{\mathcal B}<\infty,$$
where the seminorm $\|\cdot\|_{j, k; l}^{\mathcal B}$ is defined in  (\ref{qq5}).

\textit{Step 3. 
}
Any two functions $f, g\in \mathcal A$ define an element $F^{f, g}\in \mathcal B^{\mathcal A}(V\times V)$ by
\begin{equation}
\label{Ffg}
F^{f, g}(u, v):=\alpha_{Ju}(f)\alpha_v(g)\in\mathcal A, \quad\forall (u, v)\in V\times V.
\end{equation}
The \textit{deformed product} $f\times_J g$ is thus defined by the  integral (\ref{integration}) of $F^{f, g}(u, v)$ as,
\begin{equation}
\label{jproduct}
f\times_J g:=\int_V\int_V\alpha_{Ju}(f)\alpha_v(g) e(u\cdot v) \, du dv.
\end{equation}

The algebra $\mathcal A$ with its deformed product $\times_J$, together with its undeformed seminorms $\{\|\cdot\|_{j, k}\}$, defines the deformed Fr\'echet algebra $\mathcal A_J$. This  is called the \textit{deformation of the algebra $\mathcal A$ (in the direction of $J$)} as a differentiable Fr\'echet algebra. 

In the following, we obtain deformation quantizations of various  algebras of functions on EH-spaces.
We may induce a torus action $\alpha$ on the algebra $C_{b}^\infty(EH)$, or  similarly on algebras $C_0^{\infty}(EH)$ and $C^{\infty}_2(EH)$, from the torus isometric action $\sigma$ of $v\in\mathbb T^2$ on the $EH$-space  (\ref{sigmaact})  by $\alpha_v f(x)=f(\sigma_{-v}(x))$ for any $f\in C_b^{\infty}(EH)$ and  $x\in EH$.

Under the choice of the covering  $\{U_N,U_S\}$, the orbit of any point $x\in EH$ lies in the same coordinate chart as $x$.   
We  assume that 
the partition of unity $h_N$ and $h_S$ only depend on the coordinate $\theta$ so that they are invariant under the torus action $\alpha$. 

One can easily show that  the torus action $\alpha$ is isometric with respect to the family of seminorms (\ref{Q2}). We also note that each of the Fr\'echet algebras $C_b^{\infty}(EH)$ and $C_0^{\infty}(EH)$   is already smooth with respect to the action $\alpha$.
Thus, each of $C_b^{\infty}(EH)$ and $C_0^{\infty}(EH)$, with the  isometric action $\alpha$, regarded as a periodic action of $V=\mathbb R^2$,  appears exactly as the  starting point as  $(\mathcal A, \{q_m\})$ of Rieffels' deformation quantization. We can carry out step 1 to step 3  and  obtain the product $\times_J$ on the respective algebras, 
\begin{equation}
\label{Jdeform}
f\times_J g:=\int_{\mathbb R^2}\int_{\mathbb R^2}\alpha_{Ju}(f)\,\alpha_{v}(g) \,e( u\cdot v) \,du dv,
\end{equation}
where  the inner product $u\cdot v$ is the  one on $\mathbb R^2$ and $J$ is a skew-symmetric linear operator on $\mathbb R^2$. In the following we assume
$J:=\begin{pmatrix}
0&-\theta\\
\theta&0\\
\end{pmatrix},$ for some  $\theta\in\mathbb R\backslash\{0\},$
and denote $\times_J$ as $\times_{\theta}$. 

The algebra $C_b^{\infty}(EH)$ with its deformed product $\times_{\theta}$, together with its undeformed family of seminorms (\ref{Q2})
defines the deformed Fr\'echet algebra $C_b^{\infty}(EH)_{\theta}$  as the \textit{deformation quantization of  $C_b^{\infty}(EH)$}. Similarly, $C_0^{\infty}(EH)_{\theta}$ is the deformation quantization of the algebra $C_0^{\infty}(EH)$.

For the Fr\'echet algebra $C^{\infty}_2(EH)$,  the torus action $\alpha$ is  isometric with respect to the family of norms $\{\|\cdot\|_{H_{k}^2}\}_{k\in\mathbb N}$, because it is isometric with respect to the Riemannian metric. We can similarly obtain the Fr\'echet algebra $C^{\infty}_2(EH)_{\theta}$ as deformation quantization of the algebra
 $C^{\infty}_2(EH)$.
\begin{remark}
For any of the algebras in our example,  the  family of seminorms $\|\cdot\|_{j, k}$ induced from $q_m$'s as in  \textit{Step 1} is equivalent to the original family of seminorms. Indeed, the torus action is defined by the normal differentiation with respect to coordinates. 
\end{remark}

There follows some immediate observations. 
\begin{lemma}
The algebra $C^{\infty}_2(EH)_{\theta}$ is an ideal of the algebra $C_b^{\infty}(EH)_{\theta}$.
\label{H}
\end{lemma}
\begin{proof}
Let $f\in C^{\infty}_2(EH)$ and $g\in C_b^{\infty}(EH)$. Considered as elements of the algebra $C_b^{\infty}(EH)$, they define $F^{f, g}\in\mathcal B^{C_b^{\infty}(EH)}(\mathbb R^2\times\mathbb R^2)$  by (\ref{Ffg}).  
We  claim that $F^{f, g}$ lies in $\mathcal B^{C_2^{\infty}(EH)}(\mathbb R^2\times\mathbb R^2)$ so that its oscillatory integral, or product of $f\times_{\theta}g$ by definition, will be  finite in the family of seminorms on $C_2^{\infty}(EH)$ and hence   $C_2^{\infty}(EH)$-valued. In fact,
\begin{eqnarray*}
\int_{EH}|F^{f, g}(u, v)(x)|^2dVol(x)
&=&
\int_{EH}|f(Ju+x)g(v+x)|^2 dVol(x)\\
&\leq&\sup_{x\in EH}|g(x)|^2\int_{EH}|f(Ju+x)|^2 dVol(x)\\
&=&\sup_{x\in EH}|g(x)|^2\int_{EH}|f(x)|^2 dVol(x)<\infty.
\end{eqnarray*}
The last equality is by the invariance of the volume form of the integration with respect to the torus isometric action. The finiteness is because $g$ is a bounded function and $f\in C_2^{\infty}(EH)$.

Higher orders can be shown as follows. For any non-negative integer $k$,  we may expand $\nabla^k(f(Ju+x) g(v+x))$ by the Leibniz rule to a summation of terms in the form of $\nabla^l f(Ju+x)\nabla^m g(v+x)$ with $l+m=k$.  By  the assumption that  $\nabla^k f$ is $L^2$-integrable for any $k$ and $\nabla^l g$ is bounded for any $l$, each term in the summation is $L^2$-integrable. Thus $\nabla^k(f(Ju+x) g(v+x))$ is $L^2$-integrable for any $k$ and $F^{f, g}(u, v)\in C^{\infty}_{2}(EH)$ for any $(u, v)\in \mathbb R^2\times \mathbb R^2$. 
As a result, the product
$f\times_{\theta} g$ is  $C^{\infty}_2(EH)$-valued and $C^{\infty}_2(EH)_{\theta}$ is an ideal.
\end{proof}

Restriction  of the  product (\ref{Jdeform}) of the algebra $C_b^{\infty}(EH)_{\theta}$ to the algebra  $C_c^{\infty}(EH)$ gives the deformed  algebra $C_c^{\infty}(EH)_{\theta}$. We see that it is closed as an algebra as follows.  For any $f,g\in C_c^{\infty}(EH)$, the integral (\ref{Jdeform})  vanishes outside  the compact set $Orb(supp(f))\cap Orb(supp(g))$, where $Orb(U):=\{\alpha_{\mathbb T^2}(x): x\in U\subset EH\}$. Therefore, $f\times_{\theta} g$ is  of compact support and $C_c^{\infty}(EH)_{\theta}$ is thus closed. We  assign the topology of inductive limit on $C_c^{\infty}(EH)_{\theta}$ from that of $C_c^{\infty}(EH)$.   Using definitions, we have
\begin{lemma}
\label{c}
 $C_c^{\infty}(EH)_{\theta}$ is an ideal   of  the algebras $C_0^{\infty}(EH)_{\theta}$ and  $C_b^{\infty}(EH)_{\theta}$. 
 \end{lemma}
 \begin{proof}
 For  $f\in C_c^{\infty}(EH)_{\theta}$ and $g\in C_b^{\infty}(EH)_{\theta}$, the integral   (\ref{Jdeform}) vanishes outside the compact set $Orb(supp(f))$. Hence $f\times_{\theta} g$ is  $C_c^{\infty}(EH)$-valued, so that $C_c^{\infty}(EH)_{\theta}$ is an ideal of the algebras  $C_b^{\infty}(EH)_{\theta}$. The proof for the algebra $C_0^{\infty}(EH)_{\theta}$ is the same. 
 \end{proof}


The  torus action  $\alpha$ as  a compact action of an abelian group  defines a 
spectral  decomposition of a function $f$ in  the algebra $C^\infty_{b}(EH)$ or $C^{\infty}_0(EH)$,
by $$f=\sum_{s} f_s,\quad f_s(x)=\exp(i s_3\phi)\, \exp(i s_4\psi)\, h_{s}(r, \theta),$$
where  $s=(s_3, s_4)\in\mathbb Z^2$, $f_s$ satisfies  $\alpha_{v}f_s=s\, e(s\cdot v) f_s, \forall v\in \mathbb T^2$, and the series converges in the topology of uniform convergence of all derivatives. 
Under the decomposition,  the product of (\ref{Jdeform}) takes a simple form  (Chapter 2, \cite{Rieffel-1993}).
Let $f=\sum_{r}f_{r}$ and $g=\sum_{s}g_{s}$, in their respective  decompositions, be both in the algebra $C_{b}^{\infty}(EH)$ (or $C_{0}^{\infty}(EH)$), then
\begin{equation}
\label{nc}
f\times_{\theta} g
=\sum_{r, s}\sigma(r, s) f_{r}g_{s}.
\end{equation}
where $\sigma(r, s):=e(\theta (r_4 s_3 -r_3 s_4))$ and $r=(r_3, r_4), s=(s_3, s_4)\in\mathbb Z^2$. The expression (\ref{nc}) can also be restricted  to the algebra $C_c^{\infty}(EH)_{\theta}$. 
 \begin{lemma}
\label{0}
$C_0^{\infty}(EH)_{\theta}$ is an ideal of $C_b^{\infty}(EH)_{\theta}$.
\end{lemma}
\begin{proof}
For any $f\in C_0^{\infty}(EH)_{\theta}$ and $g\in C_b^{\infty}(EH)_{\theta}$, it suffices to show that $f\times_{\theta}g\in C_0^{\infty}(EH)_{\theta}$. For $g$ being zero, this is trivial.  We thus assume that $g$ is nonzero. The convergence of the series (\ref{nc}) implies that
for any $\varepsilon/2>0$, there exists an integer $N$ such that
$$|f\times_{\theta} g(x)|<\left|\sum_{|r|, |s|\leq N} \sigma(r, s) f_r(x) g_s(x)\right|+\frac{\varepsilon}{2}, $$
for any $x\in EH$, where $|r|:=|r_3|+|r_4|$ and $|s|:=|s_3|+|s_4|$. 

Since  $f_r\in C^{\infty}_0(EH)$, for each $|r|\leq N$, there exists a compact set $K(f_r)\subset EH$ such that  $$|f_r(x)|<\frac{\varepsilon}{2C}, \quad \forall x\in EH\backslash K(f_r),$$
for any  fixed constant $C$.

Therefore, for any $\varepsilon>0$, we may choose $N$ and $K(f_r)$ as above and define the union of finitely many compact sets as  $K:=\cup _{|r|\leq N} K(f_r)$, so that $x\in EH\backslash K$ implies that
$$
|f\times_{\theta} g(x)|<\left|\sum_{|r|, |s|\leq N} \sigma(r, s) f_r(x) g_s(x)\right|+\frac{\varepsilon}{2}<\frac{\varepsilon\, A_N}{2C}\sup_{x\in EH}|\sigma(r, s)g(x)| +\frac{\varepsilon}{2}, 
$$
where $A_N$ is a finite non-negative integer counting numbers of  indices $r$ and $ s$ satisfying $|r|, |s|\leq N$. If we fix the constant $C=\sup_{x\in EH}|\sigma(r, s)g(x)| A_N$, then the above inequalities give
$|f\times_{\theta} g(x)|<\varepsilon, $ whenever $x\in EH\backslash K.$
Therefore, $f\times_{\theta} g$ is  $C^{\infty}_0(EH)$-valued, and $C^{\infty}_0(EH)_{\theta}$ is an ideal of $C_b^{\infty}(EH)_{\theta}$. 
\end{proof}

We will end this subsection by introducing local algebras.
\begin{definition} \cite{rennie-2003}
 \label{localalg}
  An algebra $\mathcal A_c$ has \textit{local units} if for every finite subset of elements $\{a_i\}_{i=1}^n\subset \mathcal A_c$, there exists $\phi\in\mathcal A_c$ such that for each $i$, $\phi\, a_i=a_i\,\phi=a_i$.
  
 Let $\mathcal A$ be a Fr\'echet algebra such that  $\mathcal A_c\subset \mathcal A$ is a dense ideal with local units, then  $\mathcal A$ is called a local algebra.
 \end{definition}

 \begin{lemma}
 \label{lemA}
 The algebra $C_c^{\infty}(EH)_{\theta}$ has local units and the  algebra $C_0^{\infty}(EH)_{\theta}$ is a local algebra.
 \end{lemma}
 
 \begin{proof}
For any finite set of elements $\{f_{\beta}\}_{\beta=1}^n\subset C_c^{\infty}(EH)_{\theta}$, there exists a compact set $K$ large enough to contain the union of supports $\cup_{\beta}supp(f_{\beta})$. Let $\phi$ be a function equal to $1$ on $K$ and decaying only with respect to the $r$-variable to zero outside $K$. 
Thus defined $\phi$ satisfies $\phi=\phi_{(0, 0)}$ in the spectral decomposition
so that
 $
\phi\times_{\theta}f_{\beta}=f_{\beta}\times_{\theta}\phi=f_{\beta}$ for all $\beta.$
Thus, $(C_c^{\infty}(EH), \times_{\theta})$ is an algebra with units. 

The fact that $C^{\infty}_c(EH)$ is dense in $C^{\infty}_0(EH)$ with respect to the topology of uniform convergence of all derivatives implies that $C^{\infty}_c(EH)_{\theta}$ is dense in $C^{\infty}_0(EH)_{\theta}$, since the  family of seminorms is not deformed.
$C^{\infty}_c(EH)_{\theta}$ is an ideal in $C^{\infty}_0(EH)_{\theta}$ by Lemma \ref{c}.
Thus $C_0^{\infty}(EH)_{\theta}$ is a local algebra. 
 \end{proof}

 Lemma 3 of \cite{rennie-2003} says that there exists \textit{a local approximate unit} $\{\phi_n\}_{n\geq 1}$ for a local algebra $(\mathcal A_c\subset) \mathcal A$. In this example, we choose a family of compact sets $K_0\subset K_1\subset \dots$ in the $EH$-space, increasing in the $r$-direction. For instance,
$$K_n:=\{x\in EH: r\leq n\}, \quad \forall n\in\mathbb N.$$
Let $\{\phi_n\}_{n\in\mathbb N}$  be a family of functions with compact support $K_{n}\subset supp(\phi_n)\subset K_{n+1}$ such  that  $\phi_n$ is constant $1$ on $K_n$ and decays only with respect to $r$  to zero on $K_{n+1}$. This gives a local approximate unit. It is not hard to see that each $\phi_i$  actually commutes with functions in  the algebra $C_0^{\infty}(EH)_{\theta}$. Furthermore, the union of  the algebras $\cup_{n\in\mathbb N} [C_0^{\infty}(EH)_{\theta}]_n$, where 
 $$[C_0^{\infty}(EH)_{\theta}]_n:=\{f\in C_0^{\infty}(EH)_{\theta}: \phi_n\times_{\theta} f=f\times_{\theta} \phi_n=f\},$$
 is the algebra $C_c^{\infty}(EH)_{\theta}$. 

\label{locsec}

\subsection{Algebras of operators and deformations  of $C^*$-algebras}
\label{algop}
\begin{definition}\cite{rennie-2003}
\label{smooth2}
A $\ast$-algebra $\mathcal A$ is smooth if it is Fr\'echet and $\ast$-isomorphic to a proper dense subalgebra $i(\mathcal A)$ of a $C^*$-algebra $A$ which is stable under the holomorphic functional calculus under suitable representation.
 \end{definition}
Recall that the $q_0$ seminorm in the family  (\ref{Q2}) is  the suprenorm $\|\cdot\|_{\infty}$, which defines  $C^*$-norms on  each of the algebra  $C_b^{\infty}(EH)$ and $C_0^{\infty}(EH)$. The  $C^*$-completion of the former is the algebra $C_b(EH)$ of bounded continuous functions.  That $C_b^{\infty}(EH)$ is stable under the holomorphic functional calculus of $C_b(EH)$ implies that $C_b^{\infty}(EH)$ is a pre-$C^*$-algebra. 

The $C^*$-completion of $C_0^{\infty}(EH)$ is the algebra  $C_0(EH)$ of continuous functions vanishing at infinity. As a nonunital Banach algebra, the holomorphic functional calculus is with respect to its unitization and with respect to holomorphic functions vanishing at $0$. $C_0^{\infty}(EH)$ is stable under the holomorphic functional calculus of $C_0(EH)$ and
hence a pre-$C^*$-algebra. Similarly, the Fr\'echet algebra $C_2^{\infty}(EH)$ is also a pre-$C^*$-algebra of the $C^*$-completion $C_0(EH)$.  We see that $C_0^{\infty}(EH)$, $C_2^{\infty}(EH)$ and $C_b^{\infty}(EH)$ are smooth algebras.
 
The deformation quantizations of $C^{\infty}_0(EH)$ and $C^{\infty}_b(EH)$ obtained before  are as differentiable Fr\'echet algebras. To realize the deformed algebras as  pre-$C^*$-algebras of some deformed $C^*$-algebra, we may represent them as operators on some Hilbert space. 
  Following the construction of   \cite{Connes-landi-2001},\cite{Connes-Dubois-Violette2002}, we may obtain their representations  on the Hilbert space $\mathcal H$ of  spinors,  by using the torus isometric action. 

 Let $C^{\infty}_*(EH)_{\theta}$ stand for the algebras $C_c^{\infty}(EH)_{\theta}$, $C_0^{\infty}(EH)_{\theta}$ or $C_b^{\infty}(EH)_{\theta}$.
 The operator representation of $C^{\infty}_*(EH)_{\theta}$ on the Hilbert space $\mathcal H$ is defined by
 \begin{equation}
  \label{fn}
  L_f^{\theta}:=\sum_{r\in\mathbb Z^2}M_{f_r}V^{\theta}_r,  \end{equation}
  where $M_{f_r}$ is the normal multiplication by $f_r$ and  $V^{\theta}_r$ is a unitary operator obtained as the evaluation the unitary operator  (\ref{actiondb}) at $\tilde t_3=2\pi\theta r_4$ and $\tilde t_4=-2\pi\theta r_3$.  That is, 
  \begin{equation}
  \label{sigmarA}
  V^{\theta}_r:=e( \theta (r_4 W_3 - r_3 W_4)),\quad \quad r=(r_3, r_4)\in\mathbb Z^2.
  \end{equation}
  \begin{remark}
  \label{geom}
  Geometrically, $ V^{\theta}_r$ is the action of parallel transporting any section by $-2\pi\theta r_3$ along the $\psi$ direction followed by a parallel transporting by $2\pi\theta r_4$ along the $\phi$ direction.
  \end{remark}
  
 With the involution on $C^{\infty}_*(EH)_{\theta}$ defined by the complex conjugation of functions, we can use the property  $(f^*)_r=(f_{-r})^*$
 and $V^{\theta}_rh_s=h_sV^{\theta}_r\sigma(r, s)$ for any simple component $h_s$ from  $\sum_s h_s$, to  show
that  the representation (\ref{fn})  is a faithful $*$-representation of $C^{\infty}_*(EH)_{\theta}$. 
  
We may define the $C^*$-norm of $C^{\infty}_*(EH)_{\theta}$ by the operator norm $\|\cdot\|_{op}$ of the representation on $\mathcal H$. The series of operators (\ref{fn}) converges uniformly in the operator norm. We denote the $C^*$-completion of the algebra $C_b^\infty(EH)_{\theta}$ by $C_b(EH)_{\theta}$. It is a deformation of $C_b(EH)$ as a  $C^*$-algebra.  One can also show that $C_b^{\infty}(EH)$ is stable under the holomorphic functional calculus of $C_b(EH)_{\theta}$ and hence a pre-$C^*$-algebra. 

The $C^*$-completion $C_0(EH)_{\theta}$ of the  algebra $C^{\infty}_0(EH)_{\theta}$ defines a deformation of $C_0(EH)$ as a $C^*$-algebra. This can also be  realized as a pre-$C^*$-algebra. For similar reasons, the Fr\'echet algebra $C_2^{\infty}(EH)_{\theta}$ can  be  realized as a pre-$C^*$-algebra with the $C^*$-completion $C_0(EH)_{\theta}$.

In the commutative case, 
one can show  that $\|\cdot\|_{op}$ is bounded by the zero-th seminorm $q_0$ in the family of seminorms (\ref{semi1}). Hence the
 $C^*$-norm is weaker than the family of seminorms (\ref{semi1}). 
To see that the same holds 
in the deformed case, we note that in Rieffel's construction, the deformed Fr\'echet algebras can be represented on the space of Schwarz functions associated with a natural inner product (page 23 \cite{Rieffel-1993}) and completed to  $C^*$-algebras. Furthermore,  the correspondent $C^*$-norm is shown to be weaker than the family of seminorms defining the Fr\'echet topology (Proposition 4.10 \cite{Rieffel-1993}). We may induce a $*$-homomorphism from the $C^*$-algebra representing on $\mathcal H$ to the $C^*$-algebra representing on the space of Schwarz functions by the identity map of  functions. Since any $*$-homomorphism is norm-decreasing, we conclude that the $C^*$-norm 
 on   $C^{\infty}_*(EH)_{\theta}$  represented on  $\mathcal H$ is also weaker than the family of seminorms (\ref{semi1}) defining the topology of uniform convergence of all derivatives. 
Both of the algebras $C_0^{\infty}(EH)_{\theta}$ and $C_b^{\infty}(EH)_{\theta}$ are smooth algebras.

\subsection{Projective modules of spinor bundles}
The link between vector bundles over compact space and  projective modules is
the Serre-Swan theorem \cite{swan-1962}.  
It  is generalized  for vector bundles of finite type, of which  there exists a finite number of open sets in the open cover of the base manifold such that the bundle is trivialized on each open set \cite{vaserstein-1986}. The smooth version of the result is as follows. 
\begin{theorem}
\label{swft}
The category of complex vector bundles of finite type over $X$ for any differentiable manifold $X$ is equivalent to the category of finitely generated projective $C_b^{\infty}(X)$-modules.
\end{theorem}
\begin{remark}
There exists an alternative version of the generalized Serre-Swan theorem  \cite{rennie-2003} for vector bundles over noncompact manifolds, proved by using certain compactification of the base manifolds.
Since the simpliest one-point compactification of the Eguchi-Hanson space gives an orbifold due to the $\mathbb Z^2$-identification,  it is not straightforward to apply the construction there. 
\end{remark}
In the following, we will  use Theorem \ref{swft} to find the  projective module associated to the spinor bundle
$\mathcal S$ of the EH-space as defined  in Section \ref{EHS}.
In the coordinate charts $U_N$ and $U_S$ of the EH-space, we may choose a partition of unity $\{h_N, h_S\}$  by
\begin{equation}
\label{hn}
h_N(x):=\cos^2\frac{\theta}{2},\quad h_S( x):=\sin^2\frac{\theta}{2}, \quad x\in EH.
\end{equation}

Recall that in the unitary basis $\{f_{\alpha}\}$  of $\mathcal S_{U_N}$ and $\{f_{\beta}'\}$ of   $\mathcal S_{U_S}$, the
 transition functions  
$P
{}_{\beta}^{\alpha}$'s and $Q
{}_{\beta}^{\alpha}$'s, such that $f_{\beta}= P
{}_{\beta}^{\alpha}f'_{\alpha}$
 and 
$f'_{\beta}= Q
{}_{\beta}^{\alpha}f_{\alpha}$, are matrix entries of $P$ in (\ref{qsn}) and $Q$ in (\ref{qns}), respectively.

The idea is to extend the basis $\{f_{\alpha}\}$  on $U_N$ across the ``north pole'' $N$ and $\{f'_{\alpha}\}$ on $U_S$ across the ``south pole'' $S$ so that one can take the summation of both extended global sections to obtain a generating set of the space of smooth bounded sections of the spinor bundle $\Gamma_b^{\infty}(\mathcal S)$. 

To extend $\{f_{\alpha}\}$ across $N$, we may rescale it by the function $h_N$,
\begin{equation}
\label{Fi}
F_{\alpha}:=
\begin{cases}
f_{\alpha}  h_N & \text{ on } U_N\\
0                   & \text{ at $N$}
\end{cases},
\end{equation}
so that $F_{\alpha}$'s now decay  to zero smoothly at $N$.
Similarly, we may rescale the basis $\{f'_{\alpha}\}$ by the function $h_S$ by defining
\begin{equation}
\label{F'i}
F'_{\alpha}:=
\begin{cases}
f'_{\alpha}  h_S & \text{ on } U_S\\
0                   & \text{ at $S$}
\end{cases}.
\end{equation}
Note that on the intersection $U_N\cap U_S$, the transition function satisfies $P
{}_{\beta}^{\alpha} h_N\rightarrow 0$ whenever $h_N\rightarrow 0,$ and similarly
$Q
{}_{\beta}^{\alpha}h_S\rightarrow 0$ whenever $h_S\rightarrow 0.$
\begin{lemma}
The set of   global sections
$\{F_{\alpha}, F'_{\alpha}\}$, where $\alpha=1, \dots, 4$, are the generating set of the space of bounded smooth  sections of the spinor bundle $\Gamma^{\infty}_b(\mathcal S)$.  
\end{lemma}
\begin{proof}
The restriction 
$\{F_{\alpha}|_{U_N}\}$ where $\alpha=1, \dots, 4$ is a basis for $\mathcal S_{U_N}$. Indeed, any section $\psi\in \Gamma^{\infty}_b(\mathcal S)$ can be written as
$\psi|_{U_N}=\psi^{\alpha} f_{\alpha}= a^{\alpha} f_{\alpha} h_N= a^{\alpha} F_{\alpha}|_{U_N}, $
where $a^{\alpha}=\psi^{\alpha}/h_N.$
Similarly, the restriction $\{F'_{\alpha}|_{U_S}\}$  gives a basis for $\mathcal S_{U_S}$, since
 any section $\psi$  can be written as
$\psi|_{U_S}=\psi'^{\alpha} f'_{\alpha}= b^{\alpha} f'_{\alpha} h_S= b^{\alpha} F'_{\alpha}|_{U_S}, $
where $b^{\alpha}=\psi'^{\alpha}/h_S.$

On the intersection,
$$F_{\alpha}|_{U_N\cap U_S}=h_N\,  P
{}_{\alpha}^{\beta}\, F'_{\beta}\, h_S^{-1},\quad F'_{\alpha}|_{U_N\cap U_S}=h_S\, Q
{}_{\alpha}^{\beta}\, F_{\beta}\, h_N^{-1}.$$
Let $\{k_N, k_S\}$ be a new partition of unity such that the  $supp(k_N)\subset U_N$ and $supp(k_S)\subset U_S$. Furthermore, $k_N$  ($k_S$, respectively) is required to decay faster than $h_N$ around $N$  ($h_S$ around $S$, respectively). 
We may choose for instance,\footnote{These functions are kindly suggested by Derek Harland.} 
$$k_N(x):=\cos^2(\frac{\pi}{2}\sin^2\frac{\theta}{2}), \quad k_S(x):=\sin^2(\frac{\pi}{2}\sin^2\frac{\theta}{2}), \quad x\in EH. $$Therefore,
$a^{\alpha} k_N\rightarrow 0\text { on } U_N,$  whenever  $h_N\rightarrow 0,$
and 
$b^{\alpha} k_S\rightarrow 0\text { on } U_S,$ whenever $h_S\rightarrow 0.$
Thus, we can extend the coefficient functions $a^{\alpha}$'s and $b^{\alpha}$'s  by zero, 
$$ 
A^{\alpha}:=
\begin{cases}
a^{\alpha} k_N &\text { on }U_N\\
0           &\text{ at $N$}
\end{cases}, \quad B^{\alpha}:=
\begin{cases}
b^{\alpha} k_S &\text { on }U_S\\
0           &\text{ at $S$}
\end{cases}.$$
so that $\psi= A^{\alpha} F_{\alpha}+B^{\alpha}F'_{\alpha}.$
In fact, 
\begin{eqnarray}
\label{generate}
 A^{\alpha} F_{\alpha}+B^{\alpha}F'_{\alpha}
=
\begin{cases}
 \psi^{\alpha}  k_N f_{\alpha}  +\psi'^{\alpha} k_S f'_{\alpha}   & \text{on }U_N\cap U_S\\
 \psi'^{\alpha} k_S f'_{\alpha}   & \text{at }N\\
\psi^{\alpha} k_N  f_{\alpha} & \text{at }S
\end{cases}
=\begin{cases}
 \psi^{\alpha}   f_{\alpha}     & \text{on }U_N\\
\psi'^{\alpha} f_{\alpha}' &\text{on }U_S
\end{cases}
\end{eqnarray}
which is the section $\psi$ in $\Gamma^{\infty}_b(\mathcal S)$. 
Therefore,  $\{F_{\alpha}, F'_{\alpha}\}$ with $\alpha=1, \dots, 4$ is a generating set of $\Gamma^{\infty}_{b}(\mathcal S)$. 
\end{proof}

By  construction, we may obtain a projection in
$M_8(C^{\infty}_b(EH))$ corresponding to the spinor bundle $\mathcal S$. 
Under the standard basis of the free $C_b^{\infty}(EH)$-module $C_b^{\infty}(EH)^8$, we define the matrix, 
\begin{equation}
\label{projection}
p:=
\begin{pmatrix}
k_N\, \mathbf 1& k_N\, P
\\
k_S\,  Q
 &k_S \, \mathbf 1
\end{pmatrix}
\end{equation}
where 
$P$ and $Q$ are $4\times 4$ complex matrices  from (\ref{qsn}) and (\ref{qns}) and $\mathbf 1$ is the four by four identity matrix. 

\begin{proposition}
$\Gamma^{\infty}_b(\mathcal S)$ is a   finitely generated projective $C^{\infty}_b(EH)$-module,
\begin{equation}
\label{swb}
C^{\infty}_b(EH)^8\, p\cong \Gamma^{\infty}_b(\mathcal S).
\end{equation}
\end{proposition}
\begin{proof}
It is easy to check that $p^2=p$ and $p=p^*$. To show that (\ref{swb}) is an   isomorphism,  any section can be represented as an element in $C_b^{\infty}(EH)^8\, p$  by  construction. Conversely,
the matrix $p$ maps any element $(t_1, \dots, t_4,   t'_1,\dots, t'_4)$ of $ C^{\infty}_b(EH)^8$ to 
$$
(
(t_1+ P
{}_{1}^{\beta}t'_{\beta})\,k_N, 
(t_2+ P
{}_{2}^{\beta}t'_{\beta})\,k_N, 
(t_3+ P
{}_{3}^{\beta}t'_{\beta})\,k_N, 
(t_4+P
{}_{4}^{\beta}t'_{\beta})\,k_N,$$
$$
(t'_1+ Q
{}_{1}^{\beta}t_{\beta})\,k_S, 
(t'_2+Q
{}_{2}^{\beta}t_{\beta})\,k_S, 
(t'_3+Q
{}_{3}^{\beta}t_{\beta})\,k_S, 
(t'_4+ Q
{}_{4}^{\beta} t_{\beta})\,k_S)
.$$
Let  $A^{\alpha}=(t^{\alpha}+ P
{}_{{\alpha}}^{\beta} t'_{\beta})k_N$ and $B^{\alpha}= (t'^{\alpha}+Q
{}_{\alpha}^{\beta}t_{\beta})k_S$, for $\alpha=1,\cdots, 4$, then the image gives a section in $\Gamma^{\infty}_b(\mathcal S)$ in the form of (\ref{generate}). Therefore, (\ref{swb}) is an isomorphism.
\end{proof}

Columns of the matrix $p=(p_{\beta}^{\alpha})$ give a generating set of $\Gamma^{\infty}_b(\mathcal S)$. We may define
$\mathcal P^{k}=(p_{1}^{k}, \cdots, p_{8}^{k})^t $ for $k=1, \dots, 8$, then any element $\psi\in C^{\infty}_b(EH)^8\, p$ can be written as
$\psi= \psi_{k} \, \mathcal P^{k}$
for  functions $\psi_{k}\in C^{\infty}_b(EH)$.

\subsection{Smooth modules}
In addition to the description of a vector bundle as a  finitely generated projective module, the integrability conditions of the sections become vital when the base manifold  is noncompact. The notion of \textit{smooth module}  \cite{rennie-2003}  is  proposed to integrate the  two aspects. We will give the relevant background from the reference. 

Let $\mathcal A_0$ be an ideal in a smooth unital algebra $\mathcal A_b$. Suppose that $\mathcal A_0$ is further a local algebra containing a dense subalgebra of local units $\mathcal A_c$. Assuming  the topology on $\mathcal A_0$ is the one making it local and the topology on $\mathcal A_b$ is the one making it smooth, if the inclusion $i:\mathcal A_0\hookrightarrow \mathcal A_b$ is continuous, then $\mathcal A_0$ is a \textit{local ideal}. It is further called \textit{essential} if $\mathcal A_0\, b=\{0\}$ for some $b\in\mathcal A_b$ implies $b=0$.

Let $\mathcal A_0$ be a closed  essential local ideal in a  smooth unital algebra $\mathcal A_b$ and $p\in M_n(\mathcal A_b)$ be a projection. By pulling back the projective modules $\mathcal E_b$ defined by $\mathcal A_b^n p$ through inclusion maps $i:\mathcal A_c\hookrightarrow\mathcal A_b$, one can define the $\mathcal A_b$-finite projective $\mathcal A_c$-module $\mathcal E_c$ by $\mathcal A_c^n p$. 
 Similarly, one can define the   $\mathcal A_b$-finite projective $\mathcal A_0$-module  $\mathcal E_0$ by $\mathcal A_0^n p$. 
 
By using the Hermitian form on the projective modules $(\xi, \eta):=\sum\xi_{k}^*\eta_{k}$, one may obtain the topology on $\mathcal E_c$ induced from  the topology of inductive limit on $\mathcal A_c$, the  Fr\'echet topology on  $\mathcal E_0$ induced from  the Fr\'echet topology on $\mathcal A_0$ and the Fr\'echet topology on $\mathcal E_b$ induced from the Fr\'echet topology on  $\mathcal A_b$. Hence one has the following continuous inclusions of projective modules,
$\mathcal E_c\hookrightarrow\mathcal E_0\hookrightarrow \mathcal E_b. $

\begin{definition}
A smooth $\mathcal A_b$-module $\mathcal E_2$ is a Fr\'echet space with a continuous action of $\mathcal A_b$ such that 
$$\mathcal E_c\hookrightarrow\mathcal E_2\hookrightarrow \mathcal E_0,$$
as linear spaces, where the inclusions are all continuous.
\end{definition}

Returning to our example, we may  choose 
$\mathcal A_c$ as $C_c^{\infty}(EH)_{\theta}$, $\mathcal A_0$ as $C_0^\infty(EH)_{\theta}$, $\mathcal A_2$ as $C_2^{\infty}(EH)_{\theta}$ and $\mathcal A_b^{\infty}$ as $C_b^{\infty}(EH)_{\theta}$.

\begin{proposition}
Assuming that $C_c^{\infty}(EH)_{\theta}$ is the algebra of units, the algebras $C_c^{\infty}(EH)_{\theta}$, $C_2^{\infty}(EH)_{\theta}$ and 
$C_0^{\infty}(EH)_{\theta}$ are all essential local ideals of $C_b^{\infty}(EH)_{\theta}$ under the topology of uniform convergence of all derivatives.
\end{proposition}
\begin{proof}
$C_0^{\infty}(EH)_{\theta}$
 is  an ideal of $C_b^{\infty}(EH)_{\theta}$ by Lemma \ref{0}. Since the topology on $C_0^{\infty}(EH)_{\theta}$ and $C_b^{\infty}(EH)_{\theta}$ are both the topology of uniform convergence of all derivatives, the inclusion $C_0^{\infty}(EH)_{\theta}\hookrightarrow C_b^{\infty}(EH)_{\theta}$  is continuous. 

To show that the ideal $C_0^{\infty}(EH)_{\theta}$ is essential, we suppose that $f\in C_b^{\infty}(EH)_{\theta}$ satisfies
$g\times_{\theta} f=0$ for all $g\in C_0(EH)_{\theta}$. 
Taking  $g=1/r$,
$g\times_{\theta} f=g\times f=0$. This  implies that $f=0$, since $1/r$ is nowhere zero. Thus, $C_0^{\infty}(EH)_{\theta}$ is an essential ideal. 

$C^{\infty}_2(EH)_{\theta}$ is an ideal of $C_b^{\infty}(EH)_{\theta}$ by Lemma \ref{H}. Similar to the proof for $C_0^{\infty}(EH)_{\theta}$, $C^{\infty}_2(EH)_{\theta}$ is further an essential ideal.

$C_c^{\infty}(EH)_{\theta}$ is an ideal of $C_b^{\infty}(EH)_{\theta}$ by Lemma \ref{c}. $C_c^{\infty}(EH)_{\theta}$ carrying the topology of inductive limit is a local essential ideal, as is implied by Corollary 7 of \cite{rennie-2003} directly. 
\end{proof}

With the differential  topologies the same as their commutative restriction, there is a chain of continuous inclusions,
\begin{equation}
\label{algebra1}
C_c^{\infty}(EH)_{\theta}\hookrightarrow C_2^{\infty}(EH)_{\theta}\hookrightarrow C_0^{\infty}(EH)_{\theta}\hookrightarrow C_b^{\infty}(EH)_{\theta}.
\end{equation}

One  may define  the following  projective modules
$C_c^{\infty}(EH)_{\theta}^8\,  p$, $C_0^{\infty}(EH)_{\theta}^8\,  p$ and $ C_b^{\infty}(EH)_{\theta}^8\,  p$
by  the projection $p$ in the form of (\ref{projection}) while considered as an element in $M_8(C_b^{\infty}(EH)_{\theta})$. It is not hard to see that $p=p^*=p^2$ still holds in this case.

The family of seminorms, say $\{Q_m\}$'s, on the projective modules is  induced from the family of seminorms on the algebra, say $\{q_m\}$'s, by composing with the Hermitian form $(\cdot, \cdot)$ on the projective modules as
$Q_m(\xi):=q_m((\xi, \xi))$ for any $\xi$ in the projective module.  The  topologies on the projective modules are  defined by the induced family of seminorms. In this way, the chain of algebras  (\ref{algebra1}) induces the chain of projective modules,
$$C_c^{\infty}(EH)_{\theta}^8\,  p\hookrightarrow C_2^{\infty}(EH)_{\theta}^8 \, p\hookrightarrow C_0^{\infty}(EH)_{\theta}^8 \, p\hookrightarrow C_b^{\infty}(EH)_{\theta}^8\,  p.$$
Note that the action of $C_b^{\infty}(EH)_{\theta}$ on  
$C_2^{\infty}(EH)_{\theta}^8\,  p$ is continuous. Indeed, if a sequence of elements $\{\xi_{\beta}\}$ in $C_2^{\infty}(EH)_{\theta}^8\,  p$  satisfies  that  $Q_m(\xi_{\beta})\rightarrow 0$ as $\beta\rightarrow\infty$, then for any $f\in C_b^{\infty}(EH)_{\theta}$, $$Q_m(\xi_{\beta}\, f)=q_m((\xi_{\beta} f, \xi_{\beta} f))=q_m(f^*(\xi_{\beta}, \xi_{\beta}) f)=q_m(f^*)Q_m(\xi_{\beta}) q_m(f)\rightarrow 0,$$
where $q_m$ stands for $||\cdot||_{H_m^2}$ defined in (\ref{sob}).
Therefore, we realize $C_2^{\infty}(EH)_{\theta}^8\, p$ as a smooth module.

\section{Nonunital spectral triples and summability}
\label{localpsum}
In this section, we define  nonunital spectral triples and  consider their summability. We also consider the  regularity and measurability of the spectral triples of the isospectral deformations of EH-spaces.

Among  normed ideals in the algebra of compact operators $\mathcal K(\mathcal H)$ on a Hilbert space $\mathcal H$, 
the \textit{Dixmier trace ideal} $\mathcal L^{1, \infty}(\mathcal H)$ is the  domain of a Dixmier trace $Tr_{\omega}$, where $\omega$  is some functional on the space of bounded sequences.
 An operator $T\in \mathcal L^{1, \infty}(\mathcal H)$ is \textit{measurable} if its Dixmier trace is independent of $\omega$ and one denotes the Dixmier trace by $Tr^+(T)$. See for example \cite{Gracia-Bondia-2001}. One may define  $\Xint- T:=Tr^+(T)$ as the noncommutative integral of $T$. Apart from the Dixmier trace ideal, the \textit{generalized Schatten ideal} $\mathcal L^{p, \infty}(\mathcal H)$ for $p>1$ are the domain of operators where the $(p, \infty)$-summability are considered. 
They are related to $\mathcal L^{1, \infty}(\mathcal H)$ in a similar fashion as various Sobolev spaces  are linked. 
If the operator $T\in \mathcal L^{p, \infty}(\mathcal H)$, then $T^p\in \mathcal L^{1, \infty}(\mathcal H)$.


Rennie  (Theorem 12,  \cite{rennie-2004}) provides a measurability criterion  of operators  from \textit{local}  nonunital spectral triples. Within the locality framework, a  generalized Connes trace theorem over  commutative geodesically complete Riemannian manifold is also given (Proposition 15,  \cite{rennie-2004}). The Dixmier trace of such measurable operator agrees with  the Wodzicki residue of the operator \cite{estrada-1998}. 

Gayral and his coworkers \cite{gayral-2006-237} carry out a detailed study on  summability   of the nonunital spectral triples   from  isospectral deformations.  Their results  are  also of a local manner.

\subsection{Nonunital spectral triples and local $(p, \infty)$-summability}
\begin{definition}\cite{rennie-2003}
\label{triple}
A nonunital spectral triple $(\mathcal A, \mathcal H, \mathcal D)$ is given by
\begin{enumerate}
\item
A representation $\pi: \mathcal A\longrightarrow \mathcal B(\mathcal H)$ of a local $\ast$-algebra $\mathcal A$, containing some algebra $\mathcal A_c$ of local units as a dense ideal, on the Hilbert space $\mathcal H$. $\mathcal A$ admits a suitable unitization $\mathcal A_b$. 
\item
A self-adjoint (unbounded, densely defined) operator $\mathcal D: dom \mathcal D\longrightarrow \mathcal H$ such that $[\mathcal D, a]$ extends to a bounded operator on $\mathcal H$ for all $a\in \mathcal A_b$ and $a\, (\mathcal D-\lambda)^{-1}$ is compact for $\lambda\notin\mathbb R$ and all $a\in\mathcal A$. This is the compact resolvant condition for  nonunital triples.
\end{enumerate}
\end{definition}
We omit $\pi$ if no ambiguity arises. 
The spectral triple is \textit{even} if there exists an operator $\chi=\chi^*$ such that $\chi^2=1$, $[\chi, a]=0$ for all $a\in\mathcal A$ and $\chi\mathcal D+\mathcal D\chi=0$. Otherwise, it is \textit{odd}.

To obtain the nonunital spectral triple of the isospectral deformation of the EH-space, 
let $\mathcal A$  be the local $\ast$-algebra $C^{\infty}_0(EH)_{\theta}$ which contains the algebra of local units $C^{\infty}_c(EH)_{\theta}$ as a dense ideal. The unitization $\mathcal A_b$ is chosen as $C_b^{\infty}(EH)_{\theta}$. 
The representation $\pi$ is defined by the representation  $L_{\bullet}^{\theta}: C_b^{\infty}(EH)_{\theta}\rightarrow\mathcal B(\mathcal H)$ from (\ref{fn}). The boundedness of $L_f^{\theta}$ where $f=\sum_r f_r$ can be seen as follows, 
$$\|L_f^{\theta}\|_{op}=\left\|\sum_{r} M_{f_r}\, V^{\theta}_r\right\|_{op}
\leq\sum_{r}\|M_{f_r}\, V^{\theta}_r\|_{op}\leq\sum_{r}\|M_{f_r}\|_{op}\leq\sum_{r}\|f_r\|_{\infty}<\infty,$$
where the summations are over $\mathbb Z^2$.

Let $\mathcal D$ be extension of the Dirac operator of the spinor bundle to the Hilbert space $\mathcal H$. Since the Eguchi-Hanson space is geodesically complete, the extended operator is self-adjoint. 
We will see in the next subsection that the operator $[\mathcal D, L^{\theta}_f]$  is of degree $0$ as  a pseudodifferential operator and hence bounded.

The operator  $\chi$ is chosen to be the chirality operator defined in (\ref{chi}), such that $\chi=\chi^*$ and $\chi^2=1$.
Since $\chi$ can be realized as a fiberwise constant matrix operating on the spinor bundle, its commutativity with respect to any $L^{\theta}_f=\sum_r M_{f_r}V^{\theta}_r$ holds.
The identity $\chi\mathcal D+\mathcal D\chi=0$ is that from the commutative geometry. 

The data $(C_0^{\infty}(EH)_{\theta}, \mathcal H, \mathcal D)$ will be a nonunital spectral spectral triple once the compact resolvant condition is shown. Before that, we consider the following proposition.

\begin{proposition}
\label{locsum}
For any $f\in C_c^{\infty}(EH)_{\theta}$, 
\begin{equation}
\label{4summ} 
L^{\theta}_f(\mathcal D-\lambda)^{-1}\in\mathcal L^{4, \infty}(\mathcal H), \quad \forall\lambda\notin \mathbb R.
\end{equation} 
\end{proposition}
\begin{proof}
The proof is a straightforward generalization of Proposition 15 of \cite{rennie-2004} and references therein.

With respect to   the  local trivializations $\{U_N, U_S\}$ of the spinor bundle $\mathcal S$ coming from the stereographic projection as before, we may show the summability of the operator  (\ref{4summ})by showing the summability of the restrictions of the operator on each trivialization. Indeed, for any $f\in C_c^{\infty}(EH)_{\theta}$, the operator $L^{\theta}_f=\sum_r M_{f_r} V^{\theta}_r$ is defined by summations of  normal multiplications by $f_r$ following parallel transporting  in the $\phi$ and $\psi$ directions, so that it is well-defined when restricted on either $U_N$ or $U_S$. We may choose the partition of unity $h_N, h_S$ as (\ref{hn}) so that each function $f$ can be decomposed as
 $f=f_N+f_S$ with $f_N\in C_c^{\infty}(U_N)$ and $f_S\in C_c^{\infty}(U_S)$. 
It suffices to show that
\begin{equation}
\label{B}
L_f^{\theta}(\mathcal D-\lambda)^{-1}\in\mathcal L^{4, \infty}(L^2(\mathcal S_{U_N})),\quad\forall f\in C_c^{\infty}(U_N),
\end{equation}
and similarly for $U_S$.


For any fixed $f\in C_c^{\infty}(U_N)_{\theta}$, we can find a positive constant $R> a$ big enough, and a constant $\Theta>0$ small enough such that the compact region defined by
$$W_{R, \Theta}:=\{x\in U_N: r\leq R, \, \theta\geq \Theta\}\subset U_N,$$ 
contains the compact support of $f$. Notice that with the restricted metric from the EH-space, the region $W_{R, \Theta}$ is a compact manifold with a boundary $\partial W_{R, \Theta}$ defined by $r=R$ and $\theta=\Theta.$
 We will fix $R$ and $\Theta$ from now on, and   write $W$ instead of $W_{R, \Theta}$ and denote the restriction of the spinor bundle $\mathcal S$ on $W_{R, \Theta}$ by $\mathcal S_W$. 
 Because the integral curve starting through any point in $W$ along the $\phi$ or $\psi$ direction still lies within $W$, the action of $L_f^{\theta}$ can be restricted on sections of the subbundle $\mathcal S_W$.

To prove (\ref{B}), it suffices to prove that
$$L_f^{\theta}(\mathcal D-\lambda)^{-1}\in\mathcal L^{4, \infty}(L^2(\mathcal S_W)).$$

Let $\widetilde{W}:=W\cup_{\partial W}(-W)$  be the invertible double of the compact manifold $W$ with boundary $\partial W$, and let
the corresponding spinor bundle be $\widetilde{\mathcal S}\rightarrow \widetilde{W}$ and the corresponding Dirac operator be  $\mathcal D_I$.  Applying the Weyl's lemma \cite{Gilkey-1994} on $\widetilde{\mathcal S}\rightarrow \widetilde{W}$ as a  vector bundle over a compact manifold without boundary, we obtain
$(\mathcal D_{I}-\lambda)^{-1}\in\mathcal L^{4, \infty}(L^2(\widetilde{\mathcal S})),$
for $\lambda\notin\mathbb R$.
That is,
\begin{equation}
\label{welylemma}
\|(\mathcal D_{I}-\lambda)^{-1}\|_{4, \infty}^{\widetilde{W}\rightarrow\widetilde{W}}< \infty, \quad\forall\lambda\notin\mathbb R,
\end{equation}
where the norm is the $(4, \infty)$-Schatten norm and  we indicate the domain and image of  operators as superscript on the norms. 

As to the action of $L_f^{\theta}$, we may extend the function $f\in C_c^{\infty}(W)$ to a function $\tilde f\in C_c^{\infty}(\widetilde{W})$ by zero. 
Correspondingly, we may extend the operator 
$L_f^{\theta}:L^2(W, \mathcal S)\rightarrow L^2(W, \mathcal S)$
to 
$$L_{\tilde f}^{\theta}: L^2(\widetilde{W}, \widetilde{\mathcal S})\rightarrow L^2(\widetilde{W}, \widetilde{\mathcal S}).$$
Using the resolvant identity
$[L^{\theta}_{\tilde f}, (\mathcal D_{I}-\lambda)^{-1}]=(\mathcal D_I-\lambda)^{-1}[\mathcal D_{I}, L^{\theta}_{\tilde f}](\mathcal D_I-\lambda)^{-1},$
we have
\begin{equation}
\label{id1}
(\mathcal D_{I}-\lambda)^{-1}\, L^{\theta}_{\tilde f}
=L^{\theta}_{\tilde f}\,(\mathcal D_{I}-\lambda)^{-1}-(\mathcal D_{I}-\lambda)^{-1}(\mathcal D_I\, L^{\theta}_{\tilde f}-L^{\theta}_{\tilde f}\, \mathcal D_I)(\mathcal D_{I}-\lambda)^{-1}.
\end{equation}
By composing $L^{\theta}_{\tilde f}$ with the restriction of sections of
 $L^2(\widetilde{W}, \widetilde{\mathcal S})$ to $L^2(W, \mathcal S)$,
we obtain an operator in the same notation,  $L^{\theta}_{\tilde f}$ mapping
 from $L^2(\widetilde{W}, \widetilde{\mathcal S})$ to $L^2(W, \mathcal S)$.
 Let $\iota: W\hookrightarrow\widetilde{W}$ be the inclusion map, 
 the composition of $\iota$ with the idenity (\ref{id1}) then gives,
\begin{equation}
\label{id2}
(\mathcal D-\lambda)^{-1}\, L^{\theta}_{\tilde f}\, \iota
=L^{\theta}_{\tilde f}\,(\mathcal D_{I}-\lambda)^{-1}\, \iota
+(\mathcal D-\lambda)^{-1}(L^{\theta}_{\tilde f}\, \mathcal D_I-\mathcal D\, L^{\theta}_{\tilde f})(\mathcal D_{I}-\lambda)^{-1}\, \iota,
\end{equation}
as operators maps from
$L^2(W, \mathcal S)$ to itself.

Applying (\ref{id2}), we obtain
\begin{eqnarray*}
&&\| L^{\theta}_{ f}\, (\mathcal D-\lambda)^{-1}\|_{4, \infty}^{W\rightarrow W}\\
&=&
\|(\mathcal D-\lambda)^{-1}\, L^{\theta}_{ f}\|_{4, \infty}^{W\rightarrow W}\\
&=&\|(\mathcal D-\lambda)^{-1}\, L^{\theta}_{\tilde f}\, \iota\|_{4, \infty}^{W\rightarrow W}\\
&=&\|L^{\theta}_{\tilde f}\,(\mathcal D_{I}-\lambda)^{-1}\, \iota
+(\mathcal D-\lambda)^{-1}(L^{\theta}_{\tilde f}\, \mathcal D_I-\mathcal D\, L^{\theta}_{\tilde f})(\mathcal D_{I}-\lambda)^{-1}\, \iota \|_{4, \infty}^{W\rightarrow W}\\
&\leq&\|L^{\theta}_{\tilde f}\,(\mathcal D_{I}-\lambda)^{-1}\, \iota\|_{4, \infty}^{W\rightarrow W}
+\|(\mathcal D-\lambda)^{-1}(L^{\theta}_{\tilde f}\, \mathcal D_I-\mathcal D\, L^{\theta}_{\tilde f})(\mathcal D_{I}-\lambda)^{-1}\, \iota\|_{4, \infty}^{W\rightarrow W}.
\end{eqnarray*}

We consider the two terms in the last line separately.
Since the inclusion $\iota$ is an isometry,  the first term is bounded as
\begin{eqnarray}
\label{propA}
\|L^{\theta}_{\tilde f}\,(\mathcal D_{I}-\lambda)^{-1}\, \iota\|_{4, \infty}^{W\rightarrow W}
&\leq&\|L^{\theta}_{\tilde f}\,(\mathcal D_{I}-\lambda)^{-1}\, \|_{4, \infty}^{\widetilde{W}\rightarrow W}\nonumber\\
&\leq&\|L^{\theta}_{\tilde f}\|_{op}^{\widetilde{W}\rightarrow W}\|(\mathcal D_{I}-\lambda)^{-1}\, \|_{4, \infty}^{\widetilde{W}\rightarrow \widetilde{W}}<\infty,
\end{eqnarray}
where $\|L^{\theta}_{\tilde f}\|_{op}^{\widetilde{W}\rightarrow W}<\infty$ is because $L^{\theta}_{\tilde f}$ is the trivial extension of the bounded operator $L^{\theta}_f$ from $L^2(W, \mathcal S)$ to itself and the finiteness of  $\|(\mathcal D_{I}-\lambda)^{-1}\, \|_{4, \infty}^{\widetilde{W}\rightarrow \widetilde{W}}$ is by (\ref{welylemma}). 
The second term is bounded as
\begin{eqnarray}
\label{propB}
&&\|(\mathcal D-\lambda)^{-1}(L^{\theta}_{\tilde f}\, \mathcal D_I-\mathcal D\, L^{\theta}_{\tilde f})(\mathcal D_{I}-\lambda)^{-1}\, \iota\|_{4, \infty}^{W\rightarrow W}\nonumber\\
&\leq&\|(\mathcal D-\lambda)^{-1}(L^{\theta}_{\tilde f}\, \mathcal D_I-\mathcal D\, L^{\theta}_{\tilde f})(\mathcal D_{I}-\lambda)^{-1}\|_{4, \infty}^{\widetilde{W}\rightarrow W}\nonumber\\
&\leq&\|(\mathcal D-\lambda)^{-1}\|_{op}^{W\rightarrow W}\, \|(L^{\theta}_{\tilde f}\, \mathcal D_I-\mathcal D\, L^{\theta}_{\tilde f})\|_{op}^{\widetilde{W}\rightarrow W}\,\|(\mathcal D_{I}-\lambda)^{-1}\|_{4, \infty}^{\widetilde{W}\rightarrow \widetilde{W}}<\infty.
\end{eqnarray}
Indeed,  the finiteness of $\|(\mathcal D-\lambda)^{-1}\|_{op}^{W\rightarrow W}$ is  by the fact that  $(D-\lambda)^{-1}$ is a bounded operator on $\mathcal S\rightarrow W$ as the restriction of the bounded operator on $L^2(\mathcal S)$. For the finiteness of $\|(L^{\theta}_{\tilde f}\, \mathcal D_I-\mathcal D\, L^{\theta}_{\tilde f})\|_{op}^{\widetilde{W}\rightarrow W}$, we have 
\begin{equation*}
\|(L^{\theta}_{\tilde f}\, \mathcal D_I-\mathcal D\, L^{\theta}_{\tilde f})\|_{op}^{\widetilde{W}\rightarrow W}
=\|[\mathcal D, L_f^{\theta}]\|_{op}^{W\rightarrow W}\leq \|[\mathcal D, L_f^{\theta}]\|_{op}^{EH\rightarrow EH}<\infty,
\end{equation*}
since $\tilde f$ extends $f$ by zero and $[\mathcal D, L_f^{\theta}]$ is bounded. 
The the finiteness of $\|(\mathcal D_{I}-\lambda)^{-1}\|_{4, \infty}^{\widetilde{W}\rightarrow \widetilde{W}}$  is again by (\ref{welylemma}).

Summation of the inequalities (\ref{propA}) and (\ref{propB}) implies that
$$\| L^{\theta}_{ f}\, (\mathcal D-\lambda)^{-1}\|_{4, \infty}^{W\rightarrow W}<\infty.$$
The proof for the coordinate patch $U_S$ is the same. 
\end{proof}
As pointed out by Rennie, Proposition \ref{locsum}  implies the compact resolvant condition. 
\begin{lemma} For any $f\in C_0^{\infty}(EH)_{\theta}$,
$L_f^{\theta}\, (\mathcal D-\lambda)^{-1}\in\mathcal K(\mathcal H)$ with $\lambda\notin\mathbb R.$
\end{lemma}
\begin{proof}
Let $\{f_{\beta}\}$ is be a sequence of functions in $C_c^{\infty}(EH)_{\theta}$, which converges to the function $f\in C_0^{\infty}(EH)_{\theta}$ in  the topology of uniform convergence, then $L^{\theta}_{f_{\beta}}$ converges to $L^{\theta}_f$ in the $C^*$-operator norm, for the norm-topology is weaker than the topology of uniform convergence. This further implies that the sequence of operators 
 $\{L_{f_{\beta}}^{\theta}\, (\mathcal D-\lambda)^{-1}\}$ converges uniformly to $ L_f^{\theta}\,  (\mathcal D-\lambda)^{-1}$
 in the operator norm. The $(4, \infty)$-summability of each $L_{f_{\beta}}^{\theta}\,  (\mathcal D-\lambda)^{-1}$ by (\ref{well}) implies that they are all compact operators. As the uniform limit of a sequence of compact operators,  $L_f^{\theta}\,  (\mathcal D-\lambda)^{-1}$ is also compact.
 \end{proof}
 In summary,  the data $(C_0^{\infty}(EH)_{\theta}, \mathcal H, \mathcal D)$ of the isospectral deformationss of the Eguchi-Hanson spaces are \textit{even nonunital spectral triples}  as in Definition \ref{triple}.


\begin{definition}\cite{rennie-2004}
A (nonunital) spectral triple $(\mathcal A, \mathcal H, \mathcal D)$ is called \textit{local}, if there exists a local approximate unit $\{\phi_n\}\subset\mathcal A_c$ for $\mathcal A$ satisfying
$$\Omega_{\mathcal D}(\mathcal A_c)=\cup_n\Omega_{\mathcal D}(\mathcal A)_n,$$
where $\Omega_{\mathcal D}(\mathcal A)_n:=
\{\omega\in\Omega_{\mathcal D}(\mathcal A): \phi_n\omega=\omega\phi_n =\omega\}.$

For $p\geq 1$, the local  spectral triple  is called \textit{local $(p, \infty)$-summable} if $a\, (\mathcal D-\lambda)^{-1}\in\mathcal L^{p, \infty}(\mathcal H),$ 
 $\lambda\notin\mathbb R$, for any $a\in \mathcal A_c.$
\end{definition}
Local $(p, \infty)$-summability implies  that (Proposition 10 \cite{rennie-2004})
\begin{equation}
\label{well}
T(1+\mathcal D^2)^{-s/2}\in\mathcal L^{p/s, \infty}(\mathcal H), \quad 1\leq Re(s)\leq p,
\end{equation}
for any  $T\in\mathcal B(\mathcal H)$ such that $T\, \phi=\phi\, T=T$ for some $\phi\in\mathcal A_c$. If $Re(s)>p$, the operator is of trace class.

In considering the (local) summability of the spectral triples, we restrict ourself on   the spectral triple $(C_c^{\infty}(EH)_{\theta}, \mathcal H, \mathcal D)$. 
\begin{lemma}
The spectral triple $(C_c^{\infty}(EH)_{\theta}, \mathcal H, \mathcal D)$ is  local  $(4, \infty)$-summable.
\end{lemma}
\begin{proof}
First we show that the spectral triple is local.
We 
may  choose the local approximate unit $\{\phi_n\}$ as defined in Section \ref{locsec} so that each of $\phi_n$ remains commutative. As operators, they act only by normal multiplication $M_{\phi_n}$ on spinors. 

Define $[C_c^{\infty}(EH)_{\theta}]_n$ to be the subalgebra of $C_c^{\infty}(EH)_{\theta}$ consisting of elements $L^{\theta}_f$ such that $L^{\theta}_f\, M_{\phi_n}=M_{\phi_n}\, L^{\theta}_f=L^{\theta}_f$,
then $C_c^{\infty}(EH)_{\theta}=\cup_{n\in \mathbb N} [C_c^{\infty}(EH)_{\theta}]_n.$
Thus $$\Omega_{\mathcal D}(C_c^{\infty}(EH)_{\theta})=\Omega_{\mathcal D}(\cup_{n\in\mathbb N}[C_c^{\infty}(EH)_{\theta}]_n)
=\cup_{n\in\mathbb N}\Omega_{\mathcal D}([C_c^{\infty}(EH)_{\theta}]_n).$$ We claim that this equals to $\cup_{n\in\mathbb N}[\Omega_{\mathcal D}(C_c^{\infty}(EH)_{\theta}]_n$, where 
$$[\Omega_{\mathcal D}(C_c^{\infty}(EH)_{\theta}]_n:=\{\omega\in \Omega_{\mathcal D}(C_c^{\infty}(EH)_{\theta}): \omega\, M_{\phi_n}=M_{\phi_n}\, \omega=\omega\}.$$

By the fact that 
the orbit of the torus action of any point $x\in K_n$ remains in $K_n$, $M_{\phi_n}\, L^{\theta}_f=L^{\theta}_f\, M_{\phi_n}$ whenever $supp(f)\subset K_n$.
That the Dirac operator preserves support  implies   $$M_{\phi_n}[\mathcal D, L^{\theta}_f]=[\mathcal D, L^{\theta}_f]M_{\phi_n}=[\mathcal D, L^{\theta}_f].$$ This further gives that $\cup_{n\in\mathbb N}\Omega_{\mathcal D}([C_c^{\infty}(EH)_{\theta}]_n)\subset\cup_{n\in\mathbb N}[\Omega_{\mathcal D}(C_c^{\infty}(EH)_{\theta}]_n$. The other direction is obvious. Therefore, 
$\Omega_{\mathcal D}(C_c^{\infty}(EH)_{\theta})=\cup_{n\in\mathbb N}[\Omega_{\mathcal D}(C_c^{\infty}(EH)_{\theta}]_n, $ 
and the spectral triple is local. 

The local $(4, \infty)$-summability of the spectral triple $(C_c^{\infty}(EH)_{\theta}, \mathcal D, \mathcal H)$ is implied by
Proposition \ref{locsum}. 
\end{proof}

\subsection{Regularity of  spectral triples}
For a given spectral triple $(\mathcal A, \mathcal D, \mathcal H)$, we can define a derivation  $\delta$ on the space of linear operators on the Hilbert space $\mathcal L(\mathcal H)$ by
$$\delta(T):=[|\mathcal D|, T], \quad  T\in \mathcal L(\mathcal H).$$

A linear  operator $T$ is in the domain of the derivation  $dom\, \delta\subset\mathcal L(\mathcal H)$, if  any $\psi\in dom(|\mathcal D|)$ implies $T(\psi)\in dom(|\mathcal D|)$. 
For any positive integer $k$, 
$T$ is 
 in the domain of the $k$-th derivation $dom\, \delta^k\subset\mathcal L(\mathcal H)$, if $\delta^{k-1}(T)\in dom\, \delta$, where  $\delta^{k-1}(T)=[|\mathcal D|, [|\mathcal D|,\dots, [|\mathcal D|, T] \dots]],$
with  $k-1$ brackets.

The intersection of  domains of $\delta$ with all possible degree 
$dom^{\infty}\delta:=\cap_{k\in\mathbb N}dom\, \delta^k$ is  the \textit{smooth domain of the derivation} $\delta$. 
When $k=0$,
$dom\, \delta^0$ is simply the space of bounded operator $\mathcal B(\mathcal H)$.
Therefore, an operator $T\in dom\delta^k$ if $\delta^k(T)$ is a bounded operator.

\begin{definition}
A spectral triple $(\mathcal A, \mathcal H, \mathcal D)$ is regular if $\Omega_{\mathcal D}(\mathcal A)\subset dom^{\infty}\delta$, where $\Omega_{\mathcal D}(\mathcal A)$ is the algebra of operators  generated by $\mathcal A$ and $[\mathcal D, \mathcal A]$. 
\end{definition}
Before considering the regularity of the spectral triple, we collect some  related properties of  operators $L_f^{\theta}$ and $\mathcal D$ as pseudodifferential operators. 
The Dirac operator $\mathcal D$ on the spinor bundle $\mathcal S$ is a first order differential operator with a principal symbol,
$$\sigma^{\mathcal D}(x, \xi)=c(\xi_j dx^j).$$
where $\xi$ as a section in the cotangent bundle $T^*(EH)$ is of coordinates $(\xi_1, \dots, \xi_4)$ with respect to the basis $\{dx^i\}$, defined in the begin of Section \ref{lcc}. 
The operator $\mathcal D^2$ is a second-order differential operator with a principal symbol
\begin{equation}
\label{sigmaD2}
\sigma^{\mathcal D^2}(x, \xi)=g(\xi, \xi)\, \mathbf 1,\end{equation}
where $g$ is the induced metric tensor on the cotangent bundle from that on the tangent bundle (\ref{matrixG}).


\begin{lemma}
\label{sigML}
The principal symbol of the  pseudodifferential operator $M_f$ is 
 \begin{equation}
 \label{principalMf}
 \sigma^{M_f}(x, \xi)=M_f(x)=diag_4(f(x)),
 \end{equation}
 where $diag_r(g)$ denotes the $r\times r$ diagonal matrix of $g$ on the diagonal.
The principal symbol of the  pseudodifferential operator $L^{\theta}_f$ is 
\begin{equation}
\label{symlf}
\sigma^{L^{\theta}_f}(x, \xi)=\sum_{r=(r_3, r_4)} M_{f_r}(x)P^{\theta}(x)
\, e(\theta\, (r_3\xi_4-r_4 \xi_3) ),
\end{equation}
where the matrix-valued function $P^{\theta}(x)=P_{c_3}\circ P_{c_4}(x)$  is defined by the composition of parallel propagators along integral curves of $\partial_{\phi}$ and $\partial_{\psi}$.
\end{lemma}
\begin{proof}
Applying $M_f$ where $f=\sum_r f_r$ on the inverse Fourier transformation of a spinor $\psi$, 
$$
 M_f \left(\frac{1}{(2\pi)^4}\int_{\mathbb R^4} e^{i x\cdot \xi}\hat \psi (\xi) d\xi\right)=\frac{1}{(2\pi)^4}\int_{\mathbb R^4}diag_4(f(x)) e^{i x\cdot \xi}\hat \psi (\xi) d\xi,
$$
we see that
 $M_f$ is an order zero classical pseudodifferential operator with principal symbol (\ref{principalMf}).

From  Remark \ref{geom},  the pointwise evaluation of the operator $L^{\theta}_f$ is 
 \begin{equation}
 \label{ltheta}
 L_f^{\theta}\psi(x)=\sum_r M_{f_r}(P_{c_{3}}\circ P_{c_{3}})(\psi(x+(0, 0, -2\pi \theta r_4, 2\pi \theta r_3))),
 \end{equation}
 where  $c_{4}$ is the integral curve of the Killing field $\partial_{\psi}$ starting at $(x_1, x_2, x_3-2\pi \theta r_4, x_4+2\pi \theta r_3)$ and ending $(x_1, x_2, x_3-2\pi \theta r_4, x_4)$, and $P_{c_4}$  is assumed to be the parallel propagator with respect to the spin connection along the $c_4$. It is    evaluated at the point $(x_1, x_2, x_3-2\pi \theta r_4, x_4)$ as a four by four matrix. Similarly,
 $c_{3}$ is the integral curve  of the Killing field $\partial_{\phi}$ starting at $(x_1, x_2, x_3-2\pi \theta r_4, x_4)$ and ending at $(x_1, x_2, x_3, x_4)$.  $P_{c_3}$ is assumed to be the parallel propagator with respect to the spin connection along the $c_3$ as defined by  (\ref{propagator}). In  (\ref{ltheta}), their composition is evaluated at the point $(x_1, x_2, x_3, x_4)$ as a four by four matrix.

 Applying $L^{\theta}_f$ on the inverse Fourier transformation  of $\psi$, 
  \begin{eqnarray*}
  L_f^{\theta}\psi(x)
 &=&\frac{1}{(2\pi)^4}\int_{\mathbb R^4} \sum_r M_{f_r}P_{c_{3}} P_{c_{4}}\, \exp(i ((x+(0, 0, -2\pi \theta r_4, 2\pi \theta r_3)))\cdot \xi)\hat \psi (\xi) d\xi\\
  \end{eqnarray*}
one obtains  the symbol of $L^{\theta}_f$.  With respect to the $\xi$ variable, the complete symbol is bounded by a constant and hence is of degree $0$ and it can be chosen to be  its principal symbol, which takes the form of (\ref{symlf}).
  
  \end{proof}


\begin{proposition}
\label{smooth}
The spectral triple $(C_0^{\infty}(EH)_{\theta}, \mathcal H, \mathcal D)$ is regular.
\end{proposition}
\begin{proof}
We write  $L_f^{\theta}$ by $f$ for notational simplicity here.
As indicated in the proof of  Proposition 20 in \cite{rennie-2003},
$f, [\mathcal D, f]\in dom^{\infty}\delta$ for any $f\in\mathcal C_0^{\infty}(EH)_{\theta}$ if and only if $f, [\mathcal D, f]\in dom_{k, l\geq 0} L^k R^l,$ where
$$L(f):=(1+\mathcal D^2)^{-1/2}[\mathcal D^2, f],\quad R(f):=[\mathcal D^2, f](1+\mathcal D^2)^{-1/2},$$
for the reason that $|\mathcal D|-(1+\mathcal D^2)^{1/2}$ is bounded.  The rest of the proof is a direct generalization of the standard method in the unital case, see for instance  \cite{Gracia-Bondia-2001}. 
Denote
$ad(\mathcal D^2)^{m}(\cdot )=[\mathcal D^2, \dots, [\mathcal D^2, \cdot ]\dots],$ with $m$ brackets, so that
$$
L^k(f)=(1+\mathcal D^2)^{-k/2}ad(\mathcal D^2)^{k}(f), \quad R^l(f)=ad(\mathcal D^2)^{l}(f)(1+\mathcal D^2)^{-l/2}, $$
where $k, l\in\mathbb N.$
Their composition is
$$L^k R^l(f)=(1+\mathcal D^2)^{-k/2}ad(\mathcal D^2)^{k+l}(f)(1+\mathcal D^2)^{-l/2}.$$

The operator $ad(\mathcal D^2)(f)=[\mathcal D^2, f]$ is of order at most $1$, since  the commutator of the principal symbols (\ref{sigmaD2}) and (\ref{symlf}) vanishes.
Similarly,  the operator $ad(\mathcal D^2)^{(k+l)}(f)$ is of order at most $k+l$. This  implies that the operator $L^k R^l (f)$ is of order at most zero and hence a bounded pseudodifferential operator on $\mathcal H$. This holds for any $k$ and $l$ in $\mathbb N$. Hence $f \in dom_{k, l\geq 0} L^k R^l,$ for any $f\in C_0^{\infty}(EH_{\theta})$. 

Since $[\mathcal D, M_f]$ is a bounded operator of degree $0$ and $V^\theta(r)$ is of degree $0$, seen from (\ref{symlf}),  $[\mathcal D, L^{\theta}_f]$  is also a bounded operator of degree $0$. 
The above proof holds if $f$ is replaced by $[\mathcal D, L^{\theta}_f]$. Thus $[\mathcal D, L^{\theta}_f]\in dom_{k, l\geq 0} L^k R^l,$ for any $f\in C_0^{\infty}(EH)_{\theta}$.
Since  $L^k R^l (T)\in dom L^{0} R^{0}=\mathcal B(\mathcal H)$ for any $k, l$ where $T\in\mathcal B(\mathcal H)$ is equivalent to $T\in dom_{k, l\geq 0} L^k R^l$ for any $k, l$, we obtain $\Omega_{\mathcal D}(C_0^{\infty}(EH)_{\theta})\subset dom^{\infty}\delta$. Hence the spectral triple is regular. 
\end{proof}

\subsection{Measurability in the nonunital case}
The following is the measurability criterion of operators from  a local nonunital spectral triple \cite{rennie-2004}. 
\begin{theorem}
\label{thm12}
Let $(\mathcal A, \mathcal H, \mathcal D)$ be a regular, local $(p, \infty)$-summable spectral triple with $p\geq 1$. Suppose that $T\in\mathcal B(\mathcal H)$ such that $\psi\, T=T\,\psi=T$ for some $\psi\geq 0$ in $\mathcal A_c$. If the limit
\begin{equation}
\label{cri}
\lim_{s\rightarrow\frac{p}{2}^+}\left(s-\frac{p}{2}\right) \, Trace\left(T(1+\mathcal D^2)^{-s}\right)
\end{equation}
exists, then the operator $T(1+\mathcal D^2)^{-p/2}$ is measurable and its Dixmier trace  
 equals to the limit up to the a factor of $2/p$, 
 \begin{equation}
 \label{Dix}
Tr^+ \left(T(1+\mathcal D^2)^{-p/2}\right)=\frac{2}{p} \lim_{s\rightarrow\frac{p}{2}^+}\left(s-\frac{p}{2}\right) Trace\left(T(1+\mathcal D^2)^{-s}\right).
 \end{equation}
\end{theorem}

Implied by \cite{gayral-2006-237},
 the operators $L_f^{\theta}(1+\mathcal D^2)^{-2}$, for $f\in C_c^{\infty}(EH)_{\theta}$ from the spectral triple $(C_0^{\infty}(EH)_{\theta}, \mathcal D, \mathcal H)$ satisfies the measurability criterion (\ref{cri}) and hence the Dixmier trace can be uniquely  defined.  We  include these contents briefly for  coherence. 
\begin{lemma}
The limit 
\begin{equation}
\label{limit}
\lim_{s\rightarrow 2^+} (s-2)Trace( L_f^{\theta}(1+\mathcal D^2)^{-s}), \quad \forall f\in C_c^{\infty}(EH)_{\theta},
\end{equation}
exists and the operator $L_f^{\theta}(1+\mathcal D^2)^{-2}$ is measurable.
\label{lem4}
\end{lemma}
\begin{proof}
Since the spectral triple satisfies the local $(4, \infty)$-summability condition, (\ref{well}) implies that $L_f^{\theta}(1+\mathcal D^2)^{-s}$ for $s>2$ is of trace class and so is $M_f(1+\mathcal D^2)^{-s}$. Since both of them are of trace class, their traces agree by Corollary 3.10 of \cite{gayral-2006-237}. Thus it suffices to show that the limit
$$\lim_{s\rightarrow 2^+}(s-2)Trace( M_f(1+\mathcal D^2)^{-s}), \quad \forall f\in C_c^{\infty}(EH),$$
exists. 
We may compute the trace  of the operator by evaluating the corresponding  kernels  of operators. 
The kernel of $M_f$ is given by 
$$
K_{M_f}(x, x')=\sum M_{f_r} \delta^g_x(x'),
$$
where $\delta^g_x(x')$ is defined by requiring $\psi(x)=\int_{EH}\delta^g_x(x')\psi(x')dVol(x')$ for all $\psi\in L^2(\mathcal S)$. For $s>2$,
\begin{eqnarray*}
&&Trace( M_f\, (1+\mathcal D^2)^{-s})\\
&=&\int \int K_{M_f}(x, x')\, K_{(1+\mathcal D^2)^{-s}}(x', x)dVol(x')dVol(x)\\
&=&\int \int tr\left(\sum_{r\in\mathbb Z^2}\, diag_4(f_{r}(x))\, \delta^g_{x}(x')\right)\, K_{(1+\mathcal D^2)^{-s}}(x', x)dVol(x')dVol(x)\\
&=&4\, \int f(x)\, K_{(1+\mathcal D^2)^{-s}}(x, x)dVol(x),
\end{eqnarray*}
where $tr$ denotes the trace of a matrix.
Applying the method of  heat kernel expansion on the Laplacian transformation of the kernel as in the proof of Theorem 6.1 \cite{gayral-2006-237}, we obtain 
 that  for $s>2$,
\begin{eqnarray*}
\lim_{s\rightarrow 2^+} (s-2)Trace( M_f\, (1+\mathcal D^2)^{-s})
&=&\frac{4}{ (2\pi)^{2}}\lim_{s\rightarrow 2^+} \frac{(s-2)  \Gamma(s-2)}{\Gamma(s)} \int f(x) dVol(x)\\
&=&\frac{4}{ (2\pi)^{2}} \int f(x)dVol(x)  
<\infty.
\end{eqnarray*}
and  this equals to
$\lim_{s\rightarrow 2^+} (s-2)Trace\left( L^{\theta}_f\, (1+\mathcal D^2)^{-s}\right)$.

Since $f$ is of compact support, we can always find a function $\phi$ of value one on the compact support of $f$ and decaying to zero only with respect to the $r$ variable so that $L^{\theta}_{\phi}=M_{\phi}$ and hence $L^{\theta}_f\, M_{\phi}=M_{\phi} L^{\theta}_f =L^{\theta}_f$ holds. 
By Theorem \ref{thm12},  the operator  $L^{\theta}_f\, (1+\mathcal D^2)^{-2}$ is measurable. 
\end{proof}
(\ref{Dix}) further implies that the Dixmier trace is
\begin{equation}
\label{DixLf}
Tr^+\left( L^{\theta}_f\, (1+\mathcal D^2)^{-2}\right)
=\frac{2}{(2\pi)^{2}}\, \int_{EH} f(x)dVol(x).\end{equation}
In the reduced commutative case the operator  $M_f (1+\mathcal D^2)^{-2}$ is measurable, and $Tr^+( M_f\, (1+\mathcal D^2)^{-2})$ equals to the right hand side of (\ref{DixLf}).

The Connes trace theorem for the unital case (Theorem 7.18 \cite{Gracia-Bondia-2001}) implies  that for a  spectral triple $(\mathcal A, \mathcal H, \mathcal D)$,
\begin{equation}
\label{tracethm}
Tr^+\left(a\,(1+\mathcal D^2)^{-p/2}\right)=\frac{1}{p(2\pi)^p}\, Wres(a\, (1+\mathcal D^2)^{-p/2}),
\end{equation}
where $\mathcal D$ is the Dirac operator of some $p$-dimensional  spin manifold,  $a\, (1+\mathcal D^2)^{-p/2}$  is considered as a  elliptic pseudodifferential operator on the complex spinor  bundle $\mathcal S$ and $Wres$ is the Wodzicki residue.

Despite a full understanding of (\ref{tracethm}) in the noncommutative nonunital case,
a  Wodzicki residue computation of   $M_f^{\theta}(1+\mathcal D^2)^{1/2}$ for $f\in C_c^{\infty}(EH)$ 
shows
 \begin{eqnarray}
\label{WresMf1}
Wres\left(M_f\,( 1+\mathcal D^2)^{-2}\right)
&=&8(2\pi)^2\, \int_{EH} f(x) dVol(x), \quad f\in C_c^{\infty}(EH)
\end{eqnarray}
Comparing with (\ref{DixLf}), (\ref{tracethm}) does hold when taking $a=f$ and $p=4$. This also serves  as an example of    Proposition 15 \cite{rennie-2004} where a geodesically complete manifold is considered.

\section{Geometric conditions}
In this section, we see how the spectral triples of the isospectral deformations of the EH-spaces fit into the proposed geometric conditions  to construct  noncompact noncommutative  spin manifolds  \cite{graciabondia-2002-0204},\cite{Gayral-Bondia-Iochum-Varilly-2004}. 

For a nonunital spectral triples $(\mathcal A, \mathcal H, \mathcal D)$ as in Definition \ref{triple}, the geometric conditions are as follows.\begin{enumerate}
 
 \item[(1)]\textit{Metric dimension.}
 There is a unique non-negative integer $p$, the metric dimension, for which
 $a\, (1+\mathcal D^2)^{-1/2}$ belongs to the generalized Schatten ideal $\mathcal L^{p, \infty}(\mathcal H)$ for $a\in \mathcal A$. Moreover, $Tr^+(a\, (1+ \mathcal D^2)^{-p/2})$ is defined and not identically zero. This $p$ is even if and only if the spectral triple is even. 
 
 \item[(2)]\textit{Regularity.}
Bounded operators $a$ and $[\mathcal D, a]$, for  $a\in \mathcal A$, lie in the smooth domain of the derivation 
 $\delta=[|\mathcal D|, \cdot].$
 
 \item[(3)]\textit{Finiteness.}
The algebra $\mathcal A$ and its preferred unitization $ \mathcal A_b$ are pre-$C^*$-algebras. There exists an ideal $\mathcal A_2$ of $\mathcal A_b$, which is also a pre-$C^*$-algebra with the same $C^*$-completion as $\mathcal A$, such that the subspace of smooth vectors in $\mathcal H$
 $$
 \mathcal H^{\infty}:=\cap_{m\in\mathbb N} dom(\mathcal D^m)$$
 is an $\mathcal A_b$ finitely generated projective $\mathcal A_2$-module.

 \item[(4)]\textit{Reality.}
 There is an antiunitary operator $J$ on $\mathcal H$, such that
 $$[a, Jb^* J^{-1}]=0,$$
 for $a, b\in\mathcal A_b$. Thus $b\mapsto Jb^* J^{-1}$ is a commuting representation on $\mathcal H$ of the opposite algebra $\mathcal A_b^{\circ}$. Moreover,
 for the metric dimension $p=4$,
 $$J^2=-1, \quad J\mathcal D=\mathcal DJ, \quad J\chi=\chi J.$$
For other dimensions, we refer to the aforementioned references.
 
 \item[(5)]\textit{First order.}
 The bounded operator $[\mathcal D, a]$ commutes with the opposite algebra representation: $[[\mathcal D, a], J b^* J^{-1}]=0$ for all $a, b\in \mathcal A_b$. 
 \item[(6)]\textit{Orientation.}
 There is a Hochschild $p$-cycle $\mathbf c$ on $\mathcal A_b$, with values in $\mathcal A_b\otimes\mathcal A_b^{\circ}$. The $p$-cycle is a finite sum of terms like $(a\otimes b^{\circ})\otimes a_1\otimes\cdots\otimes a_p,$
and its natural representation $\pi_\mathcal D(\mathbf c)$  on $\mathcal H$ is defined by 
 $$\pi_\mathcal D((a_0\otimes b_0^{\circ})\otimes a_1\otimes\cdots \otimes a_p):=a_0 J b_0^* J^{-1}[\mathcal D, a_1]\cdots [\mathcal D, a_k].$$
 
 The \textit{volume form} $\pi_{\mathcal D}(\mathbf c)$  solves the equation
 $\pi_{\mathcal D}(\mathbf c)=\chi $
 in the even case and $\pi_{\mathcal D}(\mathbf c)=1$ in the odd case. 
 \end{enumerate}

\subsection{Metric dimensions}
One might show $p=4$ for the triples $(C_0^{\infty}(EH)_{\theta}, \mathcal H, \mathcal D)$ by considering the measurability of the operator $L^{\theta}_f(1+\mathcal D^2)^{-2}$ for $f\in C_0^{\infty}(EH)_{\theta}$. However, the algebra $C_0^{\infty}(EH)_{\theta}$ is not integrable, which is necessary for the computation of the Wodzicki residue \cite{estrada-1998} of the operator $L^{\theta}_f(1+\mathcal D^2)^{-2}$. Thus   $L^{\theta}_f(1+\mathcal D^2)^{-2}$ may not be measurable. Nonetheless, Lemma \ref{lem4} implies that 
operators $L^{\theta}_f(1+\mathcal D^2)^{-2}$  for $f\in C_c^{\infty}(EH)_{\theta}$ are measurable. The Dixmier trace is evaluated as
$$Tr^+(L^{\theta}_f(1+\mathcal D^2)^{-2})=\frac{2}{(2\pi)^{2}}\, \int f\,dVol,$$
which is finite and nonzero. We do not know whether this remains true for some  general  integrable algebras, for instance  $C^{\infty}_2(EH)_{\theta}$,
 lying between $C_c^{\infty}(EH)_{\theta}$ and $C_0^{\infty}(EH)_{\theta}$. 
\subsection{Finiteness}
By the construction of the ideal $C^{\infty}_2(EH)$ in Section \ref{integrable}, 
we see that  the $C_b^{\infty}(EH)$ projective $C^{\infty}_2(EH)$-module $C^{\infty}_2(EH)^8\,  p$, with $p$ as in (\ref{projection}), is the smooth domain of the Dirac operator in  $\mathcal H$. In the deformed case, we recall that $C^{\infty}_2(EH)_{\theta}^8\, p$ is  a $\mathcal C_b^{\infty}(EH)_{\theta}$ projective $C^{\infty}_2(EH)_{\theta}$-module. 

By matching  generators, we have  the isomorphism between the finitely generated  projective modules,
$C^{\infty}_2(EH)_{\theta}^8\,  p\cong C^{\infty}_2(EH)^8\,  p.$ Therefore, 
$$\cap_{m\in\mathbb N} dom(\mathcal D^m)\cong C^{\infty}_2(EH)_{\theta}^8\, p.$$
From Section \ref{algop}, the Fr\'echet algebra $C^{\infty}_2(EH)_{\theta}$ is a pre-$C^*$-algebra with the same $C^*$-completion $C_0(EH)_{\theta}$ as that of the algebra $C_0^{\infty}(EH)_{\theta}$. Hence the finiteness condition is satisfied.

As an application of  a general construction considering smooth projective modules in \cite{rennie-2003}, we may define a $\mathbb C$-valued inner product on the projective module. Since the Hermitian form 
on the  projective module $C_c^{\infty}(EH)_{\theta}^8\, p$ is $C_c^{\infty}(EH)_{\theta}$-valued,  composing with the Dixmier trace, one may define an inner product on $C_c^{\infty}(EH)_{\theta}^8\, p$ by
 $$\tau( \xi, \eta):=Tr^+\left(L^{\theta}_{(\xi|\eta)}(1+\mathcal D^2)^{-2}\right)=\frac{2}{(2\pi)^{2}} \int (\xi|\eta) \,dVol, $$
 where the equality is by (\ref{DixLf}). Here  
 the image $(\xi|\eta)=\sum \xi_k{}^*\times_{\theta}\eta_k\in C_c^{\infty}(EH)_{\theta}$ is considered as a function in $C_c^{\infty}(EH)$. 
 
 One can further take the Hilbert space completion $\overline{ C_c^{\infty}(EH)_{\theta}^8\, p}^{\tau} $  with respect to the inner product $\tau$. 
 When restricted to the commutative case,
 the inner product is simply the $L^2$-inner product on the spinor bundle, and the Hilbert space $\overline{ C_c^{\infty}(EH)^8\, p}^{\tau} $ is  the Hilbert space $\mathcal H$, appearing  in the spectral triple.

\subsection{Regularity}
 The regularity condition is implied by Proposition \ref{smooth}. 

\subsection{Reality}
The proof of the reality condition is based on the lecture notes \cite{Varilly-2006}. 
With respect to the decomposition of spinor bundle $\mathcal S=\mathcal S^+\oplus\mathcal S^-$ as in Section \ref{EHS}, we have the corresponding Hilbert space completions under the inner product coming from the $L^2$-norms, and their sum
is the Hilbert space completion of $\mathcal S$,
$
\mathcal H=\mathcal H^+\oplus\mathcal H^-.
$
Any element $\psi\in\mathcal H$ can thus be decomposed as
$\psi=(\psi^+, \psi^-)^t$. 
The  operator $J$ defined on the spinor  bundle (\ref{charge}) can be extended to the Hilbert space as an antiunitary operator  $J:\mathcal H\rightarrow \mathcal H$ by
$$J
\begin{pmatrix}
\psi^+\\
\psi^-
\end{pmatrix}
:=
\begin{pmatrix}
-\overline{\psi}^-\\
\overline{\psi}^+
\end{pmatrix},
$$
satisfying $J^2=-1$. 

We   define the representation of the opposite algebra $\mathcal A_b^{\circ}$ of $\mathcal A_b=C_b^{\infty}(EH)_{\theta}$ on $\mathcal H$, $R^{\theta}_{\bullet}: \mathcal A_{b}^{\circ}\rightarrow\mathcal B(\mathcal H)$ by 
$R^{\theta}_h:=J\, L^{\theta}_h{}^* J^{-1}.$
Specifically, for $h=\sum_s h_s$, the representation is
$$
R^{\theta}_h= \sum_s J M_{h^*_s}\, V^{\theta}_{-s} J^{-1}
=\sum_s M_{h_s}\, V^{\theta}_{-s}.
$$
 The commutativity of operators $L^{\theta}_f$ and $R^{\theta}_h$ where $f=\sum_r f_r$  is seen as follows,
 \begin{eqnarray}
 \label{form5}
 [L^{\theta}_f, R^{\theta}_h]
 &=&\sum_{r, s}f_r\, V^{\theta}_r\, h_s\, V^{\theta}_{-s}-h_s\, V^{\theta}_{-s}\, f_r\, V^{\theta}_r\nonumber\\
 &=&\sum_{r, s}f_r\, h_s \sigma(r, s)\, V^{\theta}_r\, V^{\theta}_{-s}-h_s\, \, f_r\sigma(-s, r) V^{\theta}_{-s}\, V^{\theta}_r\nonumber\\
  &=&\sum_{r, s}[f_r,  h_s ]\sigma(r, s)\, V^{\theta}_{r-s}=0,
   \end{eqnarray}
 where  identities $\sigma(r, s)=\sigma(-s, r)$ and  $V^{\theta}_rV^{\theta}_{-s}=V^{\theta}_{-s}V^{\theta}_r=V^{\theta}_{r-s}$ are applied. 
 
 As in the commutative case, 
   $\mathcal D\,  J=J\, \mathcal D$ and $J\, \chi=\chi\, J$ where  $\chi$ is the chirality operator  (\ref{chi}).

\subsection{First order}
The proof of the first order condition is again from \cite{Varilly-2006}.
For any $f=\sum_r f_r$ and $h=\sum_s h_s$ in $C^{\infty}_b(EH)_{\theta}$, the first order property  $[[\mathcal D, f_r], h_s]=0$ in the commutative case implies that,
$$
 [[\mathcal D, L^{\theta}_f], R^{\theta}_h]=\sum_{r, s}[[\mathcal D, f_r]\, V^{\theta}_r, h_s\, V^{\theta}_{-s}]
  =\sum_{r, s}[[\mathcal D, f_r], h_s] \sigma(r, s)\, V^{\theta}_{r-s}=0.
$$

\subsection{Orientation}
In Riemannian geometry, the volume form determines the orientation of a manifold. Translated to the spectral triple language, the volume form is replaced by a Hochschild cycle $\mathbf c$ which can be  represented  on $\mathcal H$ such that $\pi_D(\mathbf c)=\chi$ in the even  case. For a detailed discussion we refer to \cite{Gracia-Bondia-2001}.

We may obtain a Hochschild $4$-cycle of the spectral triple from  the classical volume form of the Eguchi-Hanson space. We will only give the construction on the coordinate chart $U_N$, that for the other chart $U_S$ is similar and the global construction can be obtained by a partition of unity. We will consider the commutative case first and then the deformed case. 

Define a new set of coordinates   by 
$u_1=x_1$, $u_2=x_2$, $u_3=e^{ ix_3}$, $u_4=e^{ i x_4}$, so that
the transition of differential forms
 $dx^i=v_j^i du^j$ is given by the diagonal matrix $V=(v_j^i):=diag(1, 1, -\frac{i}{u_3}, -\frac{i}{u_4})$. 
 Composing with the   $\vartheta^{\alpha}=h^{\alpha}_i dx^i$ where $h^{\alpha}_i $ are components of   the matrix $H$  in (\ref{matrixH}),
the transition of differential forms
$\vartheta^{\alpha}=k_{i}^{\alpha} du^i$ is given by the matrix $K=(k_i^\alpha):= H V$.
 In components, 
\begin{equation}
\label{kij}
k_1^{\alpha}=h^{\alpha}_1, \quad k_2^{\alpha}=h^{\alpha}_2, \quad 
k_3^{\alpha}=h^{\alpha}_3\,  \frac{-i}{u_3}, \quad k_4^{\alpha}=h^{\alpha}_4\, \frac{-i}{u_4},\quad\alpha=1, \cdots, 4.
\end{equation}

Similarly, the transition $du^j=\tilde v_i^j dx^i$ is given by the inverse matrix $V^{-1}=(\tilde v_i^j)$ of $V$. Composing with $dx^j=\tilde h^j_{\beta} \vartheta^\beta$ where  $\tilde h^j_{\beta}$ are elements of the inverse matrix $H^{-1}$ in  (\ref{matrixH1}) , we obtain
$du^i=\tilde k_{\beta}^i\vartheta^{\beta}$ with $ \tilde k_{\beta}^i$ as the elements of the inverse matrix $K^{-1}=  V^{-1} H^{-1}.$
In components, $$\tilde k_{\beta}^1=\tilde h^1_{\beta}, \quad \tilde k_{\beta}^2=\tilde h^2_{\beta},\quad\tilde k_{\beta}^3=i\, u_3 \, \tilde h^3_{\beta},\quad \tilde k_{\beta}^4=i\, u_4 \,  \tilde h^4_{\beta}, \quad\beta=1, \cdots, 4.$$
To avoid ambiguity, if the $u$-coordinates and $x$-coordinates appear in the same formula, we will distinguish them by adding $'$ to indices of the $u$-coordinates.
By  tensor transformations, we may obtain the Dirac operator satisfying $\mathcal D(s)=-i \gamma^{j'} \nabla_{j'}^{\mathcal S}s $ in the coordinates $\{u_i'\}$'s from (\ref{Diracx})  in the coordinates $\{x_i\}$'s  as,
\begin{eqnarray*}
\mathcal D
&=&-i\, \tilde h^{1'}_{\eta}\, \gamma^{\eta}\, \left(\partial_{1'}-\frac{1}{4}\widetilde{\Gamma}_{1\alpha}^\beta\gamma^{\alpha}\gamma_{\beta}\right)
-i\,\tilde h^{2'}_{\eta}\, \gamma^{\eta}\,\left(\partial_{2'}-\frac{1}{4}\widetilde{\Gamma}_{2\alpha}^\beta\gamma^{\alpha}\gamma_{\beta}\right)\\
&&+u_{3'} \, \tilde h^{3'}_{\eta}\, \gamma^{\eta}\, \left(\partial_{3'}+\frac{1}{4}\frac{i}{u_{3'}}\, \widetilde{\Gamma}_{3\alpha}^\beta\gamma^{\alpha}\gamma_{\beta}\right)
+u_{4'} \, \tilde h^{4'}_{\eta}\, \gamma^{\eta}\,\left(\partial_{4'}+\frac{1}{4}\frac{i}{u_{4'}}\, \widetilde{\Gamma}_{4\alpha}^\beta\gamma^{\alpha}\gamma_{\beta}\right),
\end{eqnarray*}
where $\widetilde{\Gamma}_{i\alpha}^{\beta}$'s are from (\ref{Gamma}) and $\gamma_{\alpha}=\gamma^{\alpha}$'s are from (\ref{pauli}).

The volume form of the Eguchi-Hanson space can be represented in the orthonormal basis on $U_N$ as
\begin{eqnarray}
\label{theta1234}
\vartheta^1\wedge\vartheta^2\wedge\vartheta^3\wedge\vartheta^4
&=&k^{1}_{i_1} du^{i_1}\wedge k^{2}_{i_2} du^{i_2}\wedge k^{3}_{i_3} du^{i_3}\wedge k^{4}_{i_4} du^{i_4}\nonumber\\
&=&k^{4}_{i_4}\,  k^{3}_{i_3} \,   k^{2}_{i_2}\,  k^{1}_{i_1}\, du^{i_1}\wedge  du^{i_2}\wedge  du^{i_3}\wedge  du^{i_4}.
\end{eqnarray}
We may define  a Hochschild $4$-cycle $\mathbf c_0$ in $\mathbf C_4(\mathcal A_b, \mathcal A_b\otimes\mathcal A_b^{\circ}))$, with $\mathcal A_b=C_b^{\infty}(EH)$ and $\mathcal A_b^{\circ}$ as the opposite algebra of $\mathcal A_b$, by
\begin{eqnarray}
\label{comchain}
\mathbf c_0&:=&\frac{1}{4!}\sum_{\sigma\in S_4}(-1)^{|\sigma|}(k^{\sigma(4)}_{i_{\sigma(4)}}\otimes 1^{\circ} ) ( k^{\sigma(3)}_{i_{\sigma(3)}}\otimes 1^{\circ} ) (  k^{\sigma(2)}_{i_{\sigma(2)}}\otimes 1^{\circ} )  (k^{\sigma(1)}_{i_{\sigma(1)}}\otimes 1^{\circ} )\nonumber\\
&&\qquad\otimes u^{i_{\sigma(1)}}\otimes u^{i_{\sigma(2)}}\otimes u^{i_{\sigma(3)}}\otimes u^{i_{\sigma(4)}},
\end{eqnarray}
where $\sigma$ is an element  in the permutation group $S_4$ and $(-1)^{|\sigma|}$ indicates the sign of the permutation. 
On the  $\mathcal A_b$-bimodule $ \mathcal A_b\otimes\mathcal A_b^{\circ}$, $\mathcal A_b$ acts as 
$a'(a\otimes b^0)a'':=a'aa''\otimes b^{\circ}, $ for $a\otimes b^{\circ}\in  \mathcal A_b\otimes\mathcal A_b^{\circ}$ and $a', a''\in\mathcal A_b. $

\begin{lemma}
\label{Hochb}
The Hochschild $4$-chain (\ref{comchain}) is a Hochschild cycle. That is,  
$b(\mathbf c_0)=0,$
where $b$ is the boundary operator of a Hochschild chain. 
\end{lemma}
\begin{proof}
Recall that the Hochschild  boundary operator $b$ acts on a  simple $n$-chain $a=(a_0\otimes b_0^{\circ})\otimes a_1\otimes\dots\otimes a_n$ in $\mathbf C_n(\mathcal A_b, \mathcal A_b\otimes\mathcal A_b^{\circ})$  by
\begin{eqnarray}
\label{bHoch}
b(a)&=&(a_0\otimes b_0^{\circ})a_1\otimes a_2\otimes\dots\otimes a_n\nonumber\\
&&
+\sum_{j=1}^{n-1}(-1)^j (a_0\otimes b_0^{\circ})\otimes a_1\otimes\dots\otimes a_j \, a_{j+1}\otimes\dots\otimes a_n\nonumber\\
&&+(-1)^n a_n(a_0\otimes b_0^{\circ})\otimes a_1\otimes\dots \otimes a_{n-1}.
\end{eqnarray}
Elements of $b(\mathbf c_0)$ are of three types.

The first type corresponds to  the second line in (\ref{bHoch}),
\begin{eqnarray*}
&& (-1)^{|\sigma|} (-1)^j(k^{\sigma(4)}_{i_{\sigma(4)}}\otimes 1^{\circ} ) ( k^{\sigma(3)}_{i_{\sigma(3)}}\otimes 1^{\circ} ) (  k^{\sigma(2)}_{i_{\sigma(2)}}\otimes 1^{\circ} )  (k^{\sigma(1)}_{i_{\sigma(1)}}\otimes 1^{\circ} )\\
&&\otimes u^{i_{\sigma(1)}}\otimes \dots\otimes u^{i_{\sigma(j)}}\, u^{i_{\sigma(j+1)}}\otimes\dots\otimes u^{i_{\sigma(4)}}.
\end{eqnarray*}
In the  summation of all $\sigma\in S_4$, each such term can be cancelled by a term from another $\sigma'$ which obtain from the composition of $\sigma$ by a transition between $\sigma(j)$ and $\sigma(j+1)$, as
\begin{eqnarray*}
&& (-1)^{|\sigma'|}(-1)^j (k^{\sigma'(4)}_{i_{\sigma'(4)}}\otimes 1^{\circ} ) ( k^{\sigma'(3)}_{i_{\sigma'(3)}}\otimes 1^{\circ} ) (  k^{\sigma'(2)}_{i_{\sigma'(2)}}\otimes 1^{\circ} ) (k^{\sigma'(1)}_{i_{\sigma'(1)}}\otimes 1^{\circ} )\\
&&\otimes u^{i_{\sigma(1)}}\otimes \dots\otimes u^{i_{\sigma(j+1)}}\, u^{i_{\sigma(j)}}\otimes\dots\otimes u^{i_{\sigma(4)}}.
\end{eqnarray*}
Indeed,
since $(-1)^{|\sigma|}=-(-1)^{|\sigma'|}$ and the elements in the first term from the bimodule are commuting,  the summation of such pairs is
\begin{eqnarray*}
&& (-1)^{|\sigma|} (-1)^j(k^{\sigma(4)}_{i_{\sigma(4)}}\otimes 1^{\circ} ) ( k^{\sigma(3)}_{i_{\sigma(3)}}\otimes 1^{\circ} ) (  k^{\sigma(2)}_{i_{\sigma(2)}}\otimes 1^{\circ} ) (k^{\sigma(1)}_{i_{\sigma(1)}}\otimes 1^{\circ} )\\
&&\otimes u^{i_{\sigma(1)}}\otimes \dots\otimes (u^{i_{\sigma(j)}}\,u^{i_{\sigma(j+1)}}-u^{i_{\sigma(j+1)}}\, u^{i_{\sigma(j)}})\otimes\dots\otimes u^{i_{\sigma(4)}}=0.
\end{eqnarray*}
It vanishes since $u_{i_{\sigma(j)}}\,u_{i_{\sigma(j+1)}}=u_{i_{\sigma(j+1)}}\, u_{i_{\sigma(j)}}$ as elements in $\mathcal A_b$. 

The second type corresponds to the first line in (\ref{bHoch}). After the $\mathcal A_b$-bimodule action from the right, it is in the following form,
\begin{eqnarray*}
\left(\left(k^{\sigma(4)}_{i_{\sigma(4)}}  k^{\sigma(3)}_{i_{\sigma(3)}} k^{\sigma(2)}_{i_{\sigma(2)}}  k^{\sigma(1)}_{i_{\sigma(1)}} u^{i_{\sigma(1)}}\right)\otimes 1^{\circ}  \right)\otimes u^{i_{\sigma(2)}}\otimes u^{i_{\sigma(3)}}\otimes u^{i_{\sigma(4)}}.
\end{eqnarray*}
The third type of component corresponds to the third line in (\ref{bHoch}). After the $\mathcal A_b$-bimodule action from the left, it is in the following form,
\begin{eqnarray*}
\left(\left(u^{i_{\sigma'(4)}} k^{\sigma'(4)}_{i_{\sigma'(1)}}  k^{\sigma'(3)}_{i_{\sigma'(3)}}   k^{\sigma'(2)}_{i_{\sigma'(2)}}  k^{\sigma'(1)}_{i_{\sigma'(1)}}\right)\otimes 1^{\circ} \right)\otimes  u^{i_{\sigma'(1)}}\otimes u^{i_{\sigma'(2)}}\otimes u^{i_{\sigma'(3)}}.\end{eqnarray*}
By commutativity of $\mathcal A_b$, the summation of all $\sigma$ of the second type and the third type cancel exactly when the permutation $\sigma'$   differs from $\sigma$ by a transition between $(\sigma(1), \sigma(2), \sigma(3), \sigma(4))$ to $(\sigma(4), \sigma(1), \sigma(2), \sigma(3)).$
Indeed,  such $\sigma$ and $\sigma'$ are of opposite sign. 
Therefore, all three types cancel in the summation of $\sigma\in S_4$,  and  $b(\mathbf c_{0})=0$. This shows that $\mathbf c_{0}$ is a Hochschild $4$-cycle. 
\end{proof}

We  define the representation $\pi_{\mathcal D}$ of the Hochschild cycle $\mathbf c_0$ on the Hilbert space by
$\pi_{\mathcal D}(a_0\otimes b_0^{\circ}\otimes a_1\otimes\dots\otimes a_4):=M_{a_0}M_{b_0}[\mathcal D, M_{a_1}][\mathcal D, M_{a_2}][\mathcal D, M_{a_3}][\mathcal D, M_{a_4}].$

\begin{proposition}
\label{orien}
The operator $\pi_{\mathcal D}(\mathbf c^0)=\chi$.
\end{proposition}
\begin{proof}
\begin{eqnarray*}
4!\,\pi_{\mathcal D}(\mathbf c_0)
&=&\sum_{\sigma\in S_4}(-1)^{|\sigma|}
M_{k^{\sigma(4)}_{i_{\sigma(4)}}} M_{ k^{\sigma(3)}_{i_{\sigma(3)}}} )M_{ k^{\sigma(2)}_{i_{\sigma(2)}}} M_{  k^{\sigma(1)}_{i_{\sigma(1)}}}\\
&&\qquad
c(d  u^{i_{\sigma(1)}})\, c(d  u^{i_{\sigma(2)}})\, c(d u^{i_{\sigma(3)}})
c(d u^{i_{\sigma(4)}})\\
&=&\sum_{\sigma\in S_4}(-1)^{|\sigma|}
M_{k^{\sigma(4)}_{i_{\sigma(4)}}} M_{ k^{\sigma(3)}_{i_{\sigma(3)}}} )M_{ k^{\sigma(2)}_{i_{\sigma(2)}}} M_{  k^{\sigma(1)}_{i_{\sigma(1)}}}\\
&&\qquad
\tilde k^{i_{\sigma(1)}}_{\alpha_1}\gamma^{\alpha_1}\, \tilde k^{i_{\sigma(2)}}_{\alpha_2}\gamma^{\alpha_2}\, \tilde k^{i_{\sigma(3)}}_{\alpha_3}\gamma^{\alpha_3}\, \tilde k^{i_{\sigma(4)}}_{\alpha_4}\gamma^{\alpha_4}\\
&=&\sum_{\sigma\in S_4}(-1)^{|\sigma|}\delta^{\sigma(1)}_{\alpha_1}\delta^{\sigma(2)}_{\alpha_2}\delta^{\sigma(3)}_{\alpha_3}\delta^{\sigma(4)}_{\alpha_4}\gamma^{\alpha_1}\gamma^{\alpha_2}\gamma^{\alpha_3}\gamma^{\alpha_4}\\
&=&\sum_{\sigma\in S_4}(-1)^{|\sigma|}\gamma^{\sigma(1)}\gamma^{\sigma(2)}\gamma^{\sigma(3)}\gamma^{\sigma(4)}=4!\,\gamma^{1}\gamma^{2}\gamma^{3}\gamma^{4}.
\end{eqnarray*}
Thus $\pi_{\mathcal D}(\mathbf c_0)=\chi$. 
\end{proof}

Now we consider the noncommutative case. 
Let $\mathcal A_{b, \theta}$ be $C^{\infty}_b(EH)_{\theta}$ and $\mathcal A_{b, \theta}^{\circ}$ be the opposite algebra. On the $\mathcal A_{b, \theta}$-bimodule $\mathcal A_{b, \theta}\otimes\mathcal A_{b, \theta}^{\circ}$, $\mathcal A_{b, \theta}$ acts as 
$a'(a\otimes b^0)a'':=(a'\times_{\theta}a\times_{\theta}a'')\otimes b^{\circ},$ for $a\otimes b^{\circ}\in  \mathcal A_{b,\theta}\otimes\mathcal A_{b, \theta}^{\circ}$ and $ a', a''\in\mathcal A_{b, \theta}. $

The Hochschild $4$-chain in $\mathbf C_4(\mathcal A_{b, \theta}, \mathcal A_{b, \theta}\otimes\mathcal A_{b, \theta}^{\circ})$ is defined by  
 \begin{equation}
\label{ncchain}
\mathbf c:=\frac{1}{4!}\sum_{\sigma\in S_4}(-1)^{|\sigma|}K^{\sigma(4)}_{i_{\sigma(4)}}\,  K^{\sigma(3)}_{i_{\sigma(3)}}\,   K^{\sigma(2)}_{i_{\sigma(2)}}\, K^{\sigma(1)}_{i_{\sigma(1)}} 
\otimes u^{i_{\sigma(1)}}\otimes u^{i_{\sigma(2)}}\otimes u^{i_{\sigma(3)}}\otimes u^{i_{\sigma(4)}},
\end{equation}
where $K_i^j$ is the corresponding element of $k_i^j$ 
 in the bimodule $\mathcal A_{b, \theta}\otimes\mathcal A_{b, \theta}^{\circ}$. They are chosen as,
$$
K_1^{4}:=\Delta(u_1)^{-1/2}\otimes 1^{\circ} , \quad K_2^1:=-\left(\frac{u_1}{2}\otimes 1^{\circ}\right)\varkappa(u_3),
\quad K_2^2:=\left(\frac{u_1}{2}\otimes 1^{\circ}\right)\varrho(u_3),
$$
$$K_3^1:=\left( -\frac{u_1}{2}\sin u_2\otimes 1^{\circ}\right)\,\varrho(u_3)\left(\frac{-i}{u_3}\otimes 1^{\circ}\right),$$
$$ K_3^2:=\left( -\frac{u_1}{2}\,\sin u_2\otimes 1^{\circ}       \right)\,\varkappa(u_3) \left(\frac{-i}{u_3}\otimes 1^{\circ}\right),$$
$$K_3^3:=\left( \frac{u_1}{2}\, \Delta(u_1)^{1/2}\, \cos u_2\otimes 1^{\circ}\right)\left(\frac{-i}{u_3}\otimes 1^{\circ}\right),$$
$$ K_4^3:=\left( \frac{u_1}{2}\, \Delta(u_1)^{1/2}\otimes 1^{\circ}\right)  \left(\frac{-i}{u_4}\otimes 1^{\circ}\right),$$
where 
 $$\Delta(u_1):=1-a^4/u_1^2, \quad\varkappa(u_3):=\frac{1}{2}\left((u_3^{1/2}\otimes (u_3^{1/2})^{\circ}+\bar u_3^{1/2}\otimes (\bar u_3^{1/2})^{\circ})\right),$$
$$\varrho(u_3):=\frac{1}{2i}\left((u_3^{1/2}\otimes (u_3^{1/2})^{\circ}-\bar u_3^{1/2}\otimes (\bar u_3^{1/2})^{\circ})\right).$$

\begin{remark}
\label{commutatize}
The choices of $K_i^j$'s are based on the following observation. If
 $e^{i r\phi} \in\mathcal A_{b, \theta}$ is of spectral homogeneous degree $r$, then $e^{i\frac{r}{2}\phi}\otimes (e^{i\frac{r}{2}\phi})^{\circ}$ as an element in the $\mathcal A_{b, \theta}$-bimodule is of the bimodule action satisfying
 $$e^{is\psi}(e^{i\frac{r}{2}\phi}\otimes (e^{i\frac{r}{2}\phi})^{\circ})=(e^{i\frac{r}{2}\phi}\otimes (e^{i\frac{r}{2}\phi})^{\circ})e^{is\psi},$$ for any $ e^{is\psi}$ of homogeneous degree $s$ in the algebra $\mathcal A_{b, \theta}$.
The same holds when $\phi$ and $\psi$ swap. 
In this way, all the $u_3$ appearing in the matrix $H$ of $K=HV$ can be ``commutatized''. 
\end{remark}

\begin{lemma}
The Hochschild $4$-chain (\ref{ncchain}) is a Hochschild cycle in $\mathbf Z_4(\mathcal A_{b, \theta}, \mathcal A_{b, \theta}\otimes\mathcal A_{b, \theta}^{\circ})$. I.e.,
$b(\mathbf c)=0,$
where $b$ is the boundary operator of a Hochschild chain. 
\end{lemma}
\begin{proof}
As in the commutative case, elements of  $b(\mathbf c)$ are of three types. 
The first type is,
\begin{eqnarray}
\label{d1}
&& (-1)^{|\sigma|} (-1)^j (K^{\sigma(4)}_{i_{\sigma(4)}}\times_{\theta} K^{\sigma(3)}_{i_{\sigma(3)}} \times_{\theta} K^{\sigma(2)}_{i_{\sigma(2)}} \times_{\theta}  K^{\sigma(1)}_{i_{\sigma(1)}}\nonumber\\
&&\otimes u^{i_{\sigma(1)}}\otimes \dots\otimes u^{i_{\sigma(j)}}\times _{\theta}u^{i_{\sigma(j+1)}}\otimes\dots\otimes u^{i_{\sigma(4)}}.
\end{eqnarray}
Firstly, from Remark \ref{commutatize}, 
we may observe that the noncommutative part of any $K_i^j$ has  only contributions from  terms like $\frac{i}{u_i}\times_{\theta}\cdot$, for  $i=3, 4$.
Secondly,  any term containing the product  $\frac{-i}{u_3}\times_{\theta} \frac{-i}{u_4}$ contains the product $u_4\times_{\theta} u_3$ and their product is,
 $$\frac{-i}{u_3}\times_{\theta} \frac{-i}{u_4}\times_{\theta}u_4\times_{\theta} u_3=e^{-i\theta}\, \frac{-i}{u_3} \frac{-i}{u_4}e^{i\theta}\, u_4u_3=-1.$$
This also holds when $3$ and $4$ swap. These observations imply that  the noncommutativity factor coming from the first line of (\ref{d1}) always cancels with the noncommutativity factor coming from the second line. Therefore, it reduces to the commutative case.  By  the same matching of $\sigma$'s in the proof Lemma \ref{Hochb} for terms of the first type, summation of all the terms of first type is zero. 

The second  type is
\begin{eqnarray*}
\label{d2}
&&
\left(K^{\sigma(4)}_{i_{\sigma(4)}} \times_{\theta} K^{\sigma(3)}_{i_{\sigma(3)}} \times_{\theta} K^{\sigma(2)}_{i_{\sigma(2)}}\times_{\theta}  K^{\sigma(1)}_{i_{\sigma(1)}}  \right)u^{i_{\sigma(1)}}\otimes u^{i_{\sigma(2)}}\otimes u^{i_{\sigma(3)}}\otimes u^{i_{\sigma(4)}}.
\end{eqnarray*}
Notice that 
$K^{\sigma(1)}_{i_{\sigma(1)}}$ commutes with $u_{i_{\sigma(1)}}$. 
The third type is
\begin{eqnarray*}
\label{d3}
&&u_{i_{\sigma'(4)}}\left(K^{\sigma'(4)}_{i_{\sigma'(4)}} \times_{\theta} K^{\sigma'(3)}_{i_{\sigma'(3)}}\times_{\theta} K^{\sigma'(2)}_{i_{\sigma'(2)}}\times_{\theta} K^{\sigma'(1)}_{i_{\sigma'(1)}}  \right)
\otimes u^{i_{\sigma'(1)}}\otimes u^{i_{\sigma'(2)}}\otimes u^{i_{\sigma'(3)}}.
\end{eqnarray*}
Notice that  $u_{i_{\sigma'(4)}}$ commutes with $K^{\sigma'(4)}_{i_{\sigma'(4)}}$. 
As in the commutative case, we may pair $\sigma$ and $\sigma'$ which are related by
$\sigma'(1)=\sigma(4), \sigma'(2)=\sigma(1), \sigma'(3)=\sigma(2), \sigma'(4)=\sigma(3)$ so that they are canceled through  the summation of $\sigma$. 
Three cases altogether give us $b(\mathbf c)=0$, and hence the proof. 
\end{proof}
We represent the Hochschild cycle $\mathbf c$ on the Hilbert space $\mathcal H$ by  $$\pi_{\mathcal D}(a_0\otimes b_0^{\circ}\otimes a_1\otimes\dots\otimes a_4):=L^{\theta}_{a_0}R^{\theta}_{b_0}[\mathcal D, L_{a_1}^{\theta}][\mathcal D, L_{a_2}^{\theta}][\mathcal D, L^{\theta}_{a_3}][\mathcal D, L^{\theta}_{a_4}],$$
for $a_0\otimes b_0^{\circ}\otimes a_1\otimes\dots\otimes a_4\in \mathbf Z_4(\mathcal A_{b, \theta}, \mathcal A_{b, \theta}\otimes\mathcal A_{b, \theta}^{\circ})$. 
A straightforward fact follows,
\begin{lemma}
$\pi_{\mathcal D}\left(\varkappa(u_3)\right)
=M_{\cos\phi}$ and $ \pi_{\mathcal D}\left(\varrho(u_3)\right)=M_{\sin\phi}.$
\label{comm}
\end{lemma}

\begin{proposition}
The operator $\pi_{\mathcal D}(\mathbf c)=\chi$.
\end{proposition}
\begin{proof}
By using the commutativity between the Dirac operator and $V^{\theta}_r$,  we can write down the formula for the commutators:
$$
[\mathcal D, L^{\theta}_{u_i}]=c(du^i), \quad [\mathcal D, L^{\theta}_{u_3}]
=c(du_3)V^{\theta}_{(1, 0)}, \quad [\mathcal D, L^{\theta}_{u_4}]=c(du_4) V^{\theta}_{(0, 1)}
$$
where $i=1, 2$. 
By Lemma \ref{comm}, all the nonvanishing representation of coefficients in the bimodule of the Hochschild cycle $\mathbf c$ are
$$
\pi_{\mathcal D} (K_1^{4})
=M_{\Delta(u_1)^{-1/2}}
, \quad
\pi_{\mathcal D}(K_2^1)=-M_{\frac{u_1}{2}}M_{\cos\phi}
, \quad\pi_{\mathcal D}(K_2^2)=M_{\frac{u_1}{2}}M_{\sin\phi}
$$
$$\pi_{\mathcal D}(K_3^1)= -M_{\frac{u_1}{2}\sin u_2}M_{\sin\phi}L^{\theta}_{\frac{-i}{u_3}}
, \quad \pi_{\mathcal D}(K_3^3)=M_{ \frac{u_1}{2}\, \Delta(u_1)^{1/2}\, \cos u_2}L^{\theta}_{\frac{-i}{u_3}}
$$
$$\pi_{\mathcal D}(K_3^4)= -M_{\frac{u_1}{2}\,\sin u_2} M_{\cos\phi} L^{\theta}_{\frac{-i}{u_3}}
, \quad\pi_{\mathcal D}(K_4^3)=M_{\frac{u_1}{2}\, \Delta(u_1)^{1/2}}L^{\theta}_{\frac{-i}{u_4}}
.$$

The representation $\pi_{\mathcal D}(\mathbf c)$ is thus
\begin{eqnarray*}\pi_{\mathcal D}(\mathbf c)
&=&\frac{1}{4!}\sum_{\sigma\in S_4}(-1)^{|\sigma|}
\pi_{\mathcal D}(K^{\sigma(4)}_{i_{\sigma(4)}} )\, \pi_{\mathcal D}( K^{\sigma(3)}_{i_{\sigma(3)}})\, \pi_{\mathcal D}(  K^{\sigma(2)}_{i_{\sigma(2)}} )\, \pi_{\mathcal D} (K^{\sigma(1)}_{i_{\sigma(1)}} )\nonumber\\
&&\qquad c(d u^{i_{\sigma(1)}})V^{\theta}_{i_{\sigma(1)}}\, c(d u^{i_{\sigma(2)}})V^{\theta}_{i_{\sigma(2)}}\, c(d  u^{i_{\sigma(3)}})V^{\theta}_{i_{\sigma(3)}}\, c(d  u^{i_{\sigma(4)}})V^{\theta}_{i_{\sigma(4)}},
\end{eqnarray*}
where $V^{\theta}_{i_{\sigma(k)}}:=V^{\theta}_{deg(u_{i_{\sigma(k)}})}$. 
For any fixed component in the summation we may compare the expression of $\pi_{\mathcal D}(K_{i}^j)$ and $[\mathcal D, L^{\theta}_{u_k}]$.  The result is that whenever there is a noncommutative factor generated by some $\pi_{\mathcal D}(K_{i}^j)$ as $V^{\theta}_{deg(1/u_i)}$ there is a corresponding noncommutative factor generated by   $[\mathcal D, L^{\theta}_{u_i}]$ as $V^{\theta}_{deg(u_i)}$. Furthermore, these paired noncommutative factors cancel consistently. Thus, each component in the summation is simply the same as that in the commutative case. Applying Proposition \ref{orien}, the summation gives $\chi$ again and this completes the proof of the orientation condition,
$\pi_{\mathcal D}(\mathbf c)=\chi$. 
\end{proof}

\section{Conclusions}
We have obtained the nonunital spectral triples of the isospectral deformations of the Eguchi-Hanson spaces along  torus isometric actions and studied   analytical properties of the triple. We have also tested the proposed  geometric conditions of a noncompact noncommutative geometry on this example.

There are possible generalizations in the following directions. Firstly,  we may further consider the Poincar\'e duality   of  nonunital spectral triples \cite{rennie-2001}. Secondly, we may take the conical singularity limit of EH-spaces and consider the spectral triple of the conifold. Thirdly, we may realize the spectral triple as a complex noncommutative geometry defined by \cite{froehlich-1998-193}. Finally, we may  deform the EH-spaces, and possibly for  more general ALE-spaces, by using the hyper-K\"ahler quotient structures.

\section*{Acknowledgements}The author  thanks  Lucio Cirio, Giovanni Landi for their interests and comments,  and   Derek Harland for helpful discussions, without which the research would have taken much longer. Finally, the author  wants to thank  Adam Rennie for sharing his insights, encouragement and comments on the  draft.

The work is supported  by the Dorothy Hodgkin Scholarship from  University of Durham.
\bibliographystyle{unsrt}
\bibliography{refadhm1}

\end{document}